\documentclass[11pt]{article} 
\usepackage{amssymb} 
\usepackage[dvips]{epsfig} 
\usepackage{graphicx} 

\newtheorem{Lemma}{Lemma}[section] 
\newtheorem{Theorem}{Theorem} 
\newtheorem{Proposition}[Lemma]{Proposition} 
 
\newtheorem{Remark}[Lemma]{Remark}

\newenvironment{Proof} 
 {\begin{trivlist} \item[]{\bf Proof. }} 
 {\QED\end{trivlist}} 
 
\newenvironment{Acknowledgment} 
 {\begin{trivlist}\item[]\textbf{Acknowledgments }}{\end{trivlist}} 
 
\makeatletter\@addtoreset{figure}{section}\makeatother

\renewcommand{\theLemma}{\arabic{section}.\arabic{Lemma}} 
\renewcommand{\theRemark}{\arabic{section}.\arabic{Remark}} 
 
\makeatletter\@addtoreset{equation}{section}\makeatother 
\renewcommand{\theequation}{\thesection.\arabic{equation}} 
 
\newcommand{\C}{\mathbb{C}} 
 
\newcommand{\R}{\mathbb{R}} 
\newcommand{\Z}{\mathbb{Z}} 
\newcommand{\T}{\mathbb{T}} 
 
\def\Re{\mathop{\mathrm{Re}}} 
\def\Im{\mathop{\mathrm{Im}}}

\newcommand{\rmO}{\mathcal{O}}

\newcommand{\rme}{\mathrm{e}} 
\newcommand{\rmi}{\mathrm{i}}

\newcommand{\per}{\mathrm{per}} 
\newcommand{\loc}{\mathrm{loc}} 
\newcommand{\sign}{\mathrm{sign}} 
\newcommand{\dist}{\mathrm{dist}} 
\newcommand{\I}{\mathcal I} 
\newcommand{\dd}{\,\mathrm{d}} 
\newcommand{\EE}{\mathcal{E}} 
\newcommand{\MM}{\mathcal{M}} 
\newcommand{\HH}{\mathcal{H}} 
\newcommand{\KK}{\mathcal{K}} 
\newcommand{\TT}{\mathcal{T}} 
\newcommand{\RR}{\mathcal{R}} 
\newcommand{\NN}{\mathcal{N}} 
\newcommand{\1}{\mathbf{1}} 
 
\renewcommand{\phi}{\varphi} 
\newcommand{\half}{{\textstyle \frac12}} 
 
\def\ze{^0_{a,b,\gamma}} 
\def\abg{_{a,b,\gamma}}

\newcommand{\reff}[1]{(\ref{#1})} 
 
\def\eqdef{\buildrel\hbox{\small{def}}\over =} 
\def\build#1_#2^#3{\mathrel{ 
  \mathop{\kern 0pt#1}\limits_{#2}^{#3}}} 
\def\QED{\mbox{}\hfill$\Box$} 
 
\newdimen\texpscorrection 
\newdimen\figcenter 
\def\figurewithtex #1 #2 #3 #4 #5\cr{\null 
  {\goodbreak\figcenter=\hsize\relax 
  \advance\figcenter by -#4truecm 
  \divide\figcenter by 2 
  \begin{figure}[hbt] 
  \vskip #3truecm\noindent\hskip\figcenter 
  \includegraphics{#1}{\hskip\texpscorrection\input #2 } 
  \vskip 0.8truecm\noindent \vbox{\noindent {\footnotesize #5}} 
  \end{figure}}} 
\def\point#1 #2 #3 {\rlap{\kern #1 truecm 
\raise #2 truecm \hbox{#3}}}

\begin{document} 
 
\thispagestyle{empty} 
 
\begin{center} 
{\LARGE Stability of small periodic waves for the nonlinear\\[1ex] 
Schr{\"o}dinger equation} 
\\[8mm] 
{\large\bf Thierry Gallay} 
\\ 
Institut Fourier\\ 
Universit{\'e} de Grenoble I\\ 
B.P. 74\\ 
38402 Saint-Martin-d'H{\`e}res, France 
\\[4mm] 
{\large\bf Mariana  H\u{a}r\u{a}gu\c{s}} 
\\ 
D\'epartement de Math\'ematiques\\ 
Universit\'e de Franche-Comt\'e\\ 
16 route de Gray\\ 
25030 Besan\c con, France 
\end{center} 
 
\vspace*{1cm} 
 
\begin{abstract} 
\noindent The nonlinear Schr{\"o}dinger equation possesses three 
distinct six-parameter families of complex-valued quasi-periodic 
travelling waves, one in the defocusing case and two in the 
focusing case. All these solutions have the property that their 
modulus is a periodic function of $x-ct$ for some $c \in \R$. In 
this paper we investigate the stability of the small amplitude 
travelling waves, both in the defocusing and the focusing case. 
Our first result shows that these waves are orbitally stable 
within the class of solutions which have the same period and the 
same Floquet exponent as the original wave. Next, we consider 
general bounded perturbations and focus on spectral stability. We 
show that the small amplitude travelling waves are stable in the 
defocusing case, but unstable in the focusing case. The 
instability is of side-band type, and therefore cannot be detected 
in the periodic set-up used for the analysis of orbital stability. 
\end{abstract} 
 
\vfill 
 
\noindent {\bf Running head:} {Periodic waves in the NLS equation} 
 
\noindent {\bf Corresponding author:} Mariana H\u{a}r\u{a}gu\c{s}, 
{\tt haragus@math.univ-fcomte.fr} 
 
 
\newpage 
 
 
\section{Introduction} 
\label{s:1} 
 
We consider the one-dimensional cubic nonlinear Schr\"odinger 
equation (NLS) 
\[ 
  \rmi U_t(x,t) + U_{xx}(x,t) \pm |U(x,t)|^2 U(x,t) \,=\, 0~, 
\] 
where $x \in \R$, $t \in \R$, $U(x,t) \in \C$, and the signs $+$ 
and $-$ in the nonlinear term correspond to the focusing and the 
defocusing case, respectively. In both cases the NLS equation 
possesses {\it quasi-periodic} solutions of the general form 
\begin{equation}\label{e:qpdef} 
  U(x,t) \,=\, \rme^{\rmi (p x -\omega t)}\,V(x-ct)~, 
  \quad x \in \R~, \quad t \in \R~, 
\end{equation} 
where $p,\omega,c$ are real parameters and the wave profile $V$ is 
a \emph{complex-valued periodic} function of its argument. The aim 
of the present paper is to investigate the stability properties of 
these particular solutions, at least when the wave profile $V$ is  
small. It turns out that the discussion is very similar in both cases,  
so for simplicity we restrict our presentation to the defocusing equation 
\begin{equation}\label{e:nls} 
  \rmi U_t(x,t) + U_{xx}(x,t) - |U(x,t)|^2 U(x,t) \,=\, 0~, 
\end{equation} 
and only discuss the differences which occur in the focusing case 
at the end of the paper. 
 
A crucial role in the stability analysis is played by the various 
symmetries of the NLS equation. The most important ones for our 
purposes are the four continuous symmetries: 
\begin{itemize} 
  \item {\it phase invariance:} $U(x,t) \mapsto 
    U(x,t)\,\rme^{\rmi\phi}$, $\phi\in\R$; 
  \item {\it translation invariance:} $U(x,t) \mapsto 
    U(x+\xi,t)$, $\xi\in\R$; 
  \item {\it Galilean invariance:} $U(x,t)\mapsto 
    \rme^{-\rmi\big(\frac v2x + \frac{v^2}4t\big)} U(x+vt,t)$, 
    $v\in\R$; 
  \item {\it dilation invariance:} $U(x,t) \mapsto \lambda 
    U(\lambda x,\lambda^2t)$, $\lambda > 0$; 
\end{itemize} 
and the two discrete symmetries: 
\begin{itemize} 
  \item {\it reflection symmetry:} $U(x,t) \mapsto 
    U(-x,t)$, 
  \item {\it conjugation symmetry:} $U(x,t) \mapsto 
    \overline{U}(x,-t)$. 
\end{itemize} 
  
As is well-known, the Cauchy problem for equation \reff{e:nls} is 
globally well-posed on the whole real line in the Sobolev space 
$H^1(\R,\C)$, see e.g. \cite{CW,GV1,GV2,K2}. Alternatively, one can  
solve the NLS equation on a bounded interval $[0,L]$ with periodic  
boundary conditions, in which case an appropriate function space is  
$H^1_\per([0,L],\C)$. In both situations, we have the 
following conserved quantities: 
\begin{eqnarray*} 
  E_1(U) &=& \frac12 \int_\I |U(x,t)|^2 \dd x~, 
    \\[2mm] 
  E_2(U) &=& \frac{\rmi}{2}\int_\I {\overline U}(x,t)U_x(x,t)\dd x~, 
   \\[2mm] 
  E_3(U)&=& \int_\I \Bigl(\,\frac12 
    |U_x(x,t)|^2 + \frac14|U(x,t)|^4\Bigr)\dd x~, 
\end{eqnarray*} 
where $\I$ denotes either the whole real line or the bounded interval 
$[0,L]$. The quantities $E_1$ and $E_2$ are conserved due to the phase 
invariance and the translation invariance, respectively, whereas the 
conservation of $E_3$ originates in the fact that equation 
\reff{e:nls} is autonomous. 
 
The symmetries listed above are also useful to understand the 
structure of the set of all quasi-periodic solutions of \reff{e:nls}. 
Assume that $U(x,t)$ is a solution of \reff{e:nls} of the form 
\reff{e:qpdef}, where $V : \R \to \C$ is a bounded function.  Since 
$|U(x,t)| = |V(x-ct)|$, the translation speed $c \in \R$ is uniquely 
determined by $U$, except if the modulus $|V|$ is constant. In any 
case, using the Galilean invariance, we can transform $U(x,t)$ into 
another solution of the form \reff{e:qpdef} with $c = 0$.  Once this 
is done, the temporal frequency $\omega$ is in turn uniquely 
determined by $U(x,t)$, except in the trivial case where $V$ is 
identically zero. In view of the dilation invariance, only the sign of 
$\omega$ is important, so we can assume without loss of generality 
that $\omega \in \{-1;0;1\}$. Setting $U(x,t) = \rme^{-\rmi \omega t} 
W(x)$, we see that $W(x) = \rme^{\rmi p x} V(x)$ is a bounded solution 
of the ordinary differential equation 
\begin{equation}\label{e:snls} 
  W_{xx}(x) + \omega W(x) -|W(x)|^2 W(x) \,=\, 0~, \quad 
  x \in \R~. 
\end{equation} 
If $\omega = 0$ or $\omega = -1$, is it straightforward to verify that 
$W \equiv 0$ is the only bounded solution of \reff{e:snls}, thus we 
assume from now on that $\omega = 1$. Equation~\reff{e:snls} is 
actually the stationary Ginzburg-Landau equation and the set of its 
bounded solutions is well-known \cite{BR1,DGJ,Ga,GH2}. There are two 
kinds of solutions of \reff{e:snls} which lead to quasi-periodic 
solutions of the NLS equation of the form \reff{e:qpdef}: 
\begin{itemize} 
\item A family of \emph{periodic solutions} with constant modulus 
$W(x) = (1-p^2)^{1/2}\,\rme^{\rmi (px+\phi)}$, where $p \in 
[-1,1]$ and $\phi \in [0,2\pi]$. The corresponding solutions of 
\reff{e:nls} are called {\it plane waves}. The general form of 
these waves is 
\[ 
  U(x,t) \,=\, \rme^{\rmi (px - \omega t)}\,V~, 
\] 
where $p \in \R$, $\omega \in \R$, and $V \in \C$ satisfy the 
dispersion relation $\omega = p^2 + |V|^2$. 
 
\item A family of \emph{quasi-periodic solutions} of the form 
$W(x) = r(x)\,\rme^{\rmi \phi(x)}$, where the modulus $r(x)$ and 
the derivative of the phase $\phi(x)$ are periodic with the same 
period. Any such solution can be written in the equivalent form 
$W(x) = \rme^{\rmi p x}\,Q(2kx)$, where $p \in \R$, $k > 0$, and 
$Q : \R \to \C$ is $2\pi$-periodic. In particular, 
\begin{equation}\label{e:WQ} 
  U(x,t) = \rme^{-\rmi t} W(x) = \rme^{\rmi (px -t)} Q(2kx) 
\end{equation} 
is a quasi-periodic solution of \reff{e:nls} of the form 
\reff{e:qpdef} (with $c = 0$ and $\omega = 1$). We shall refer to 
such a solution as a {\it periodic wave}, because its profile 
$|U(x,t)|$ is a (non-trivial) periodic function of the space 
variable $x$. Important quantities related to the periodic wave 
\reff{e:WQ} are the period of the modulus $T = \pi/k$, and the 
Floquet multiplier $\rme^{\rmi pT}$. For small amplitude solutions 
($|Q| \ll 1$) the minimal period $T$ is close to $\pi$, hence $k 
\approx 1$, and the Floquet multiplier is close to $-1$, so that 
we can choose $p \approx 1$. 
\end{itemize} 
 
While the plane waves form a three-parameter family, we will see 
in Section~\ref{s:ex} that the periodic waves form a six-parameter 
family of solutions of \reff{e:nls}. However, the number of 
independent parameters can be substantially reduced if we use the 
four continuous symmetries listed above. Indeed it is easy to 
verify that any plane wave is equivalent either to $U_1(x,t) = 0$ 
or to $U_2(x,t) = \rme^{-\rmi t}$. In a similar way, the set 
of all periodic waves reduces to a two-parameter family. 
 
As far as the stability of the {\it plane waves} is concerned, the 
conserved quantities $E_1$ and $E_3$ immediately show that the 
trivial solution $U_1 = 0$ is stable (in the sense of Lyapunov) with 
respect to perturbations in $H^1(\R)$ or $H^1_\per([0,L])$, for 
any $L > 0$. The same conservation laws also imply that the plane 
wave $U_2 = \rme^{-\rmi t}$ is {\it orbitally stable} in the following sense. 
Assume that $\I = [0,L]$ is a bounded interval, and let $U(x,t)$ 
be the solution of \reff{e:nls} with initial data $U(x,0) = 1 + 
V_0(x)$, where $V_0 \in H^1_\per(\I)$ and $\|V_0\|_{H^1(\I)} \le 
\epsilon$. If $\epsilon > 0$ is small enough, then 
\begin{equation}\label{e:planestab} 
  \inf_{\phi \in [0,2\pi]} \|U(\cdot,t) - \rme^{\rmi \phi} 
  \|_{H^1(\I)} \,\le\, C(\I)\epsilon~, \quad \hbox{for all } 
  t \in \R~, 
\end{equation} 
where the constant $C(\I)$ depends only on the length of the 
interval $\I$. This stability property is easily established 
using the conserved quantity 
\[ 
  E(U)\,=\,\int_\I \Bigl(\,\frac12 |U_x(x,t)|^2 + 
  \frac14(|U(x,t)|^2 - 1)^2\Bigr)\dd x \,=\, 
  E_3(U) - E_1(U) + \frac14 |\I|~. 
\] 
A similar result holds for small perturbations of $U_2$ in 
$H^1(\R)$. In that case, the bound \reff{e:planestab} holds for 
any bounded interval $\I \subset \R$, but the conservation of 
$E(U)$ does not prevent the norm $\|U(\cdot,t) - 
\rme^{-\rmi t}\|_{H^1(\R)}$ from growing as $|t| \to \infty$. We refer 
to \cite[Section 3.3]{Z} for a detailed analysis of the stability 
of plane waves. 
 
The stability question is much more difficult for {\it periodic 
  waves}. In contrast to dissipative systems for which nonlinear 
stability of periodic patterns has been established for rather general 
classes of perturbations, including localized ones (see e.g. 
\cite{Sch}), no such result is available so far for dispersive 
equations. In the particular case of NLS, the stability of the ground 
state solitary waves has been intensively studied (see e.g. 
\cite{CL,W}), but relatively little seems to be known about the 
corresponding question for periodic waves. A partial spectral analysis 
is carried out by Rowlands \cite{R}, who shows that the periodic waves 
are unstable in the focusing case and stable in the defocusing case, 
provided disturbances lie in the long-wave regime. Spectral stability 
has also been addressed for certain NLS-type equations with periodic 
potentials \cite{BR2,OSY}.  Very recently, Angulo \cite{A} has shown 
that the family of ``dnoidal waves'' of the focusing NLS equation is 
orbitally stable with respect to perturbations which have the same 
period as the wave itself. In all these previous works, the wave 
profile $V$ is assumed to be real-valued. Here we restrict ourselves 
to small amplitude solutions, but allow for general complex-valued 
wave profiles.  While the nonlinear stability of these waves with 
respect to bounded or localized perturbations remains an open problem, 
we treat here two particular questions: orbital stability with respect 
to periodic perturbations, and spectral stability with respect to 
bounded or localized perturbations. 
 
Our first result shows that the periodic waves of \reff{e:nls} are 
{\it orbitally stable} within the class of solutions which have 
the same period and the same Floquet multiplier as the original 
wave: 
 
\begin{Theorem}\label{th:orbit} {\rm (Orbital stability)}\\ 
Let $X = H^1_\per([0,2\pi],\C)$. There exist positive constants 
$C_0$, $\epsilon_0$, and $\delta_0$ such that the following holds. 
Assume that $W(x) = \rme^{\rmi px}\,Q_\per(2kx)$ is a solution of 
\reff{e:snls} with $Q_\per \in X$, $\|Q_\per\|_X \le \delta_0$, 
and $p,k \approx 1$. For all $R \in X$ such that $\|R\|_X \le 
\epsilon_0$, the solution $U(x,t) = \rme^{\rmi (px-t)}Q(2kx,t)$ of 
\reff{e:nls} with initial data $U(x,0) = \rme^{\rmi 
px}(Q_\per(2kx) + R(2kx))$ satisfies, for all $t \in \R$, 
\begin{equation}\label{e:orbit} 
  \inf_{\phi,\xi \in [0,2\pi]} \|Q(\cdot,t) - \rme^{\rmi \phi} 
  Q_\per(\cdot - \xi)\|_X \,\le\, C_0 \|R\|_X~. 
\end{equation} 
\end{Theorem} 
 
\medskip\noindent{\bf Remarks}\\ 
{\bf 1.} We point out that Theorem~\ref{th:orbit} holds uniformly 
for all quasi-periodic solutions of \reff{e:nls} with small 
amplitude. In particular the unperturbed solution $\rme^{\rmi 
(px-t)} \,Q_\per(2kx)$ can be either a periodic wave, or a plane 
wave, 
or even the zero solution. \\[1mm] 
{\bf 2.} The proof of Theorem~\ref{th:orbit} relies on the 
classical approach to orbital stability which goes back to 
Benjamin \cite{Be} (see also \cite{ABS,Bo,W}). While for solitary 
waves this method gives a rather complete answer to the stability 
question, in the case of periodic waves it allows to prove orbital 
stability only if we restrict ourselves to solutions which have 
the same periodicity properties as the original wave (see however 
Remark~\ref{r:op} below for a discussion of this limitation). In 
this paper we shall use the general framework developed by 
Grillakis, Shatah, and Strauss \cite{GSS1,GSS2}, with appropriate 
modifications to obtain a uniform stability result for small 
waves. Note that a direct application of the stability theorem in 
\cite{GSS2} would give the same conclusion as in 
Theorem~\ref{th:orbit}, but with stability constants $C_0$ and 
$\epsilon_0$ depending on the wave 
profile $Q_\per$.\\[1mm] 
{\bf 3.} Following the approach of \cite{GSS2} it is shown in 
\cite{GH2} that all periodic waves of \reff{e:nls} are orbitally 
stable in the sense of \reff{e:orbit}, without any restriction 
on the amplitude of the wave profile $Q_\per$. The argument in 
\cite{GH2} relies in part on the results obtained in the present 
paper, and uses a global parametrization of the set of quasi-periodic 
solutions of \reff{e:snls} which is very different from the 
explicit series expansions that we use here to describe the small 
amplitude solutions.\\[1mm] 
{\bf 4.} It is worth considering what Theorem~\ref{th:orbit} 
exactly means in the particular case where $W$ is a {\it  
real-valued} periodic solution of \reff{e:snls} (such a solution 
is often referred to as a ``cnoidal wave'' in the literature).  
In that case we have $W(x) = \rme^{\rmi px}\,Q_\per(2kx)$  
where $p = k = \pi/T$ and $T > \pi$ is the minimal period  
of $|W|$. The Floquet multiplier $\rme^{\rmi pT}$ is therefore  
exactly equal to $-1$, so that $W(x+T) = -W(x)$ for all  
$x \in \R$. In particular, we see that $W$ is periodic with  
(minimal) period $L = 2T$. Thus Theorem~\ref{th:orbit} shows  
that the $L$-periodic cnoidal wave $U(x,t) = \rme^{-\rmi t}W(x)$ 
is orbitally stable with respect to $L$-periodic perturbations 
$\widetilde W$ {\it provided} that $\widetilde W(x+L/2) =  
-\widetilde W(x)$ for all $x \in \R$. As is explained in \cite{A},  
without this additional assumption the classical approach does not  
allow to prove the stability of cnoidal waves with respect to periodic  
perturbations.  
 
\medskip 
Next, we investigate the spectral stability of the periodic waves 
with respect to bounded, or localized, perturbations. Although 
spectral stability is weaker than nonlinear stability, it provides 
valuable information about the the linearization of the system at 
the periodic wave. Our second result is: 
 
\begin{Theorem}\label{th:spec} {\rm (Spectral stability)}\\ 
Let $Y = L^2(\R,\C)$ or $Y=C_b(\R,\C)$. There exists $\delta_1 > 
0$ such that the following holds. Assume that $W(x) = \rme^{\rmi 
px}\,Q_\per(2kx)$ is a solution of \reff{e:snls} with $Q_\per \in 
X$, $\|Q_\per\|_X \le \delta_1$, and $p,k \approx 1$, just as in 
Theorem~\ref{th:orbit}. Then the spectrum of the linearization of 
\reff{e:nls} about the periodic wave $\rme^{-\rmi t}W(x)$ in the 
space $Y$ entirely lies on the imaginary axis. Consequently, this 
wave is spectrally stable in $Y$. 
\end{Theorem} 
 
The proof of Theorem~\ref{th:spec} is based on the so-called 
Bloch-wave decomposition, which reduces the spectral study of the 
linearized operator in the space $Y$ to the study of the spectra of a 
family of linear operators in a space of periodic functions. Bloch 
waves are well-known for Schr{\"o}dinger operators with periodic 
potentials \cite{RS} and have been extensively used in dissipative 
problems \cite{Mi,Sc1, Sc2, Sch}, but also in a number of dispersive 
problems \cite{BR2,HLS,OSY}. The advantage of such a decomposition is 
that the resulting operators have compact resolvent, and therefore 
only point spectra. The main step in the analysis consists in locating 
these point spectra. For our problem, we rely on perturbation 
arguments for linear operators in which we regard the operators 
resulting from the Bloch-wave decomposition as small perturbations of 
operators with constant coefficients. The latter ones are actually 
obtained from the linearization of \reff{e:nls} about zero, and 
Fourier analysis allows to compute their spectra explicitly.  The 
restriction to small amplitudes is essential in this perturbation 
argument, and we do not know whether spectral stability holds for 
large waves. 
 
The rest of the paper is organized as follows. In 
Section~\ref{s:ex}, we briefly describe the set of all bounded 
solutions of \reff{e:snls}, and we introduce an analytic 
parametrization of the small amplitude solutions which will be 
used throughout the paper. In Section~\ref{s:os} we recall the 
main ideas of the orbital stability method, and we apply it with 
appropriate modifications to prove Theorem~\ref{th:orbit}. 
Spectral stability is established in Section~\ref{s:ss}, using 
Bloch-wave decomposition and the pertubation argument mentioned 
above. 
 
Finally, in Section~\ref{s:foc} we discuss the stability of the small 
periodic waves of the focusing NLS equation. In contrast to the 
defocusing case, the focusing NLS equation possesses two different 
families of quasi-periodic solutions of the form \reff{e:qpdef}, one 
for $\omega>0$ and the other for $\omega<0$ \cite{GH2}. Small 
solutions exist only within the first family, and their stability 
properties can be analyzed as in the defocusing case. However, while 
for periodic perturbations we obtain the same orbital stability result 
as in Theorem~\ref{th:orbit}, it turns out that the small periodic 
waves are spectrally {\it unstable} in the focusing case. Unstable 
spectrum is detected for perturbations with wave-numbers which are 
close to that of the original wave (side-band instability). 
As for the second family, which contains only large waves, we 
refer to \cite{A,GH2} for a proof of orbital stability and 
to \cite{OSY} for a discussion of spectral stability.  
 
\begin{Acknowledgment} 
The authors thank A. De Bouard and L. Di Menza for fruitful 
discussions. This work was partially supported by the French 
Ministry of Research through grant ACI JC~1039. 
\end{Acknowledgment}

\section{Parametrization of small periodic waves} 
\label{s:ex} 
 
In this section, we briefly review the bounded solutions of 
equation \reff{e:snls} with $\omega=1$: 
\begin{equation}\label{e:sgl} 
  W_{xx}(x) + W(x) -|W(x)|^2 W(x) \,=\, 0~,\quad W : \R \to \C~, 
\end{equation} 
and we give a convenient parametrization of all small solutions. 
If we interpret the spatial variable $x \in \R$ as a ``time'', 
equation \reff{e:sgl} becomes an integrable Hamiltonian dynamical 
system with two degrees of freedom. The conserved quantities 
are the ``angular momentum'' $J$ and the ``energy'' $E$: 
\begin{equation}\label{e:JEdef} 
  J \,=\, \Im(\overline{W}W_x)~, \qquad 
  E \,=\, \frac12|W_x|^2 + \frac12|W|^2 - \frac14|W|^4~. 
\end{equation} 
If $W$ is a solution of \reff{e:sgl} with $J \neq 0$, then $W(x) 
\neq 0$ for all $x \in \R$, so that we can introduce the polar 
coordinates $W(x) = r(x)\,\rme^{\rmi \phi(x)}$. The invariants 
then become 
\[ 
  J \,=\, r^2 \phi_x~, \qquad E \,=\,\frac{r_x^2}2 + 
  \frac{J^2}{2r^2} + \frac{r^2}2 - \frac{r^4}4~. 
\] 
 
The set of bounded solutions of \reff{e:sgl} can be entirely 
described in terms of these two invariants \cite{BR1,DGJ,Ga,GH2}. 
In the parameter space $(J,E)$ there is an open set 
\begin{equation}\label{e:Ddef} 
  D \,=\, \Bigl\{(J,E)\in\R^2 \,\Big|\, J^2 < 4/27\,,~E_-(J) < 
  E < E_+(J)\Bigr\}~, 
\end{equation} 
where $E_-,E_+$ are explicit functions of $J$, such that the 
closure $\overline D$ consists of all values of $(J,E)$ which give 
rise to bounded solutions $W$ of \reff{e:sgl} (see Fig.~1). 
Furthermore, we have the following classification for $(J,E)$  
in $\overline D$: 
\begin{enumerate} 
\item If $E = E_-(J)$, then $W$ is a periodic solution with 
constant modulus and linear phase, i.e. $W(x) = W_{p,\phi}(x) = 
(1-p^2)^{1/2}\,\rme^{\rmi (px+\phi)}$ with $1/3 \le p^2 \le 1$ and 
$\phi \in [0,2\pi]$. 
 
\item If $E = E_+(J)$, then either $W = W_{p,\phi}$ for some $p^2 
\le 1/3$ and some $\phi \in [0,2\pi]$, or $W$ is a homoclinic 
orbit connecting $W_{p,\phi_-}$ at $x = -\infty$ to $W_{p,\phi_+}$ 
at $x = +\infty$ for some $\phi_-,\phi_+ \in [0,2\pi]$. 
 
\item If $E_-(J) < E < E_+(J)$ and $J \neq 0$, then the modulus 
and the phase derivative of $W$ are both periodic with the same 
period $T(J,E) > \pi$. Let $\Phi(J,E)$ be the increment of the 
phase over a period of the modulus, so that $W(x+T) = 
\rme^{\rmi \Phi}W(x)$ for all $x \in \R$. In general $\Phi$ is 
not a rational multiple of $\pi$, hence the solution $W$ of 
\reff{e:sgl} is typically not periodic, but only quasi-periodic. 
In the particular case where $J=0$, then $\rme^{\rmi \Phi} = -1$ 
and $W$ is periodic with period $2T(0,E)$. 
\end{enumerate} 
 
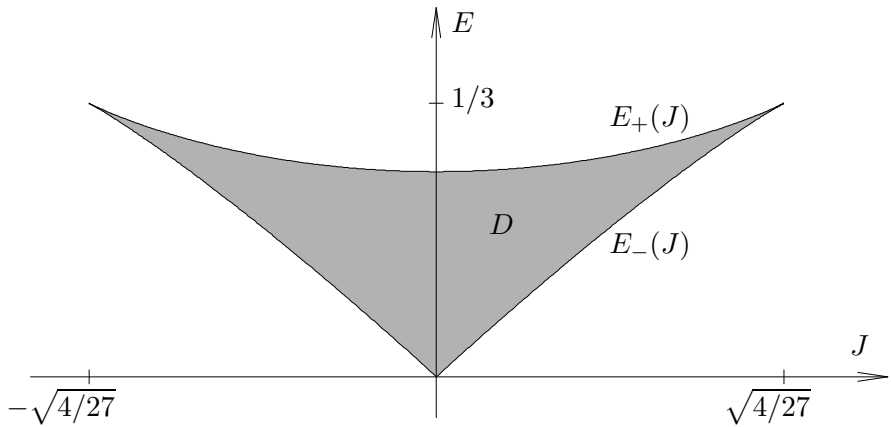
\figurewithtex fig1.ps Fig1.tex 6.000 12.000 {\bf Fig.~1:} The 
region $D$ in the parameter space $(J,E)$ for which \reff{e:sgl} 
has bounded solutions.\cr 
 
\noindent For a fixed pair $(J,E)\in \overline D$, the bounded 
solution $W$ of \reff{e:sgl} satisfying \reff{e:JEdef} is unique 
up to a translation and a phase factor. In case (iii), the period 
$T$ and the phase increment  $\Phi$ (or the Floquet multiplier 
$\rme^{\rmi \Phi}$) are important quantities which play a crucial 
role in the stability analysis of the quasi-periodic solutions of 
\reff{e:sgl}, both for the Schr\"odinger and the Ginzburg-Landau 
dynamics. A number of properties of $T$ and $\Phi$ are collected 
in \cite{GH2}. In particular, if we define the renormalized phase 
\begin{equation}\label{e:Psidef} 
  \Psi(J,E) \,=\, \left\{ 
   \begin{array}{ccl} \Phi(J,E) - \pi\,\sign(J) & \hbox{if} 
     & J \neq 0~,\\ 0 & \hbox{if} & J = 0~, 
   \end{array} 
  \right. 
\end{equation} 
then $T : D\to\R$ and $\Psi : D \to \R$ are smooth functions 
of $(J,E) \in D$, in contrast to $\Phi(J,E)$ which is discontinuous 
at $J = 0$. In addition, $T \approx \pi$ and $\Psi \approx 0$ 
for small solutions $W\approx0$. 
 
The periodic solutions $W_{p,\phi}$ of \reff{e:sgl} correspond to 
plane waves of the NLS equation. We are mainly interested here in 
the quasi-periodic solutions described in (iii) above, which 
correspond to periodic waves of the form \reff{e:qpdef}. To see 
this, fix $(J,E) \in D$ and let $W : \R \to \C$ be a bounded 
solution of \reff{e:sgl} satisfying \reff{e:JEdef}. We set 
\begin{equation}\label{e:Prep} 
  W(x) \,=\, \rme^{\rmi \ell x}\,P(kx)~, \quad x \in \R~, 
\end{equation} 
in which $k$ and $\ell$ are related to the period $T(J,E)$ and the 
renormalized phase $\Psi(J,E)$ through 
\begin{equation}\label{e:kldef} 
  k \,=\, \frac{\pi}{T(J,E)}~, \quad \hbox{and}\quad 
  \ell \,=\, \frac{\Psi(J,E)}{T(J,E)}~. 
\end{equation} 
As $|W(x)| = |P(kx)|$ is $T$-periodic (in $x$) by the definition 
of $T(J,E)$, it is clear that $|P(y)|$ is $\pi$-periodic (in $y$). 
Moreover, since $W(x+T) = \rme^{\rmi \Phi}W(x) = -\rme^{\rmi 
\Psi}W(x)$, we also have $P(y+\pi) = -P(y)$ for all $y \in \R$, 
hence $P$ is $2\pi$-periodic. Thus $U(x,t) = \rme^{-\rmi t}W(x) = 
\rme^{\rmi \ell x}\rme^{-\rmi t}P(kx)$ is a quasi-periodic 
solution of \reff{e:nls} of the form \reff{e:qpdef}, with $\omega 
= 1$ and $c = 0$. 
 
\begin{Remark} 
Using the continuous symmetries of the NLS equation, we obtain 
for each pair $(J,E) \in D$ a four-parameter family of 
periodic waves: 
\[ 
  U_{c,\lambda,\phi,\xi}(x,t) \,=\, \lambda \,\rme^{\rmi p_{c,\lambda}x} 
  \,\rme^{-\rmi \omega_{c,\lambda}t}\,\rme^{\rmi \phi} 
  P(k\lambda(x-ct)-\xi)~, 
\] 
where $c \in \R$, $\lambda > 0$, $\phi,\xi \in [0,2\pi]$, and 
$p_{c,\lambda} = \lambda\ell + c/2$, $\omega_{c,\lambda} = 
\lambda^2 + c\lambda\ell +c^2/4$. Taking into account the 
parameters $J,E$, we obtain altogether a six-parameter family of 
periodic waves of \reff{e:nls}. 
\end{Remark} 
 
Alternatively, we can write the solution \reff{e:Prep} of 
\reff{e:sgl} in the form 
\begin{equation}\label{e:Qrep} 
  W(x) \,=\, \rme^{\rmi (\ell + k)x}\,Q^+(2kx) \,=\, \rme^{\rmi 
  (\ell - k)x}\,Q^-(2kx)~, \quad x \in \R~, 
\end{equation} 
where $Q^\pm(z) = \rme^{\mp \rmi z/2}P(z/2)$. By construction, 
$Q^\pm$ and $|Q^\pm|$ are now periodic functions with the 
{\it same} minimal period $2\pi$. The representation \reff{e:Qrep} 
turns out to be more convenient than \reff{e:Prep} to study the 
orbital stability of the periodic waves in the next section. 
 
The global parametrization of the quasiperiodic solutions of 
\reff{e:sgl} in terms of the invariants $(J,E)$ is natural, but it 
is not very convenient as far as small solutions are concerned 
because the admissible domain $D$ is not smooth near the origin 
(see Fig.~1). For this reason, we now introduce an analytic 
parametrization of the {\it small} solutions of \reff{e:sgl}. We 
start from the representation \reff{e:Prep}, and we choose as 
parameters the first nonzero Fourier coefficients of the 
$2\pi$-periodic function $P$: 
\[ 
  a \,=\, \frac{1}{2\pi}\int_0^{2\pi} P(y)\,\rme^{\rmi y}\dd y~, \quad 
  b \,=\, \frac{1}{2\pi}\int_0^{2\pi} P(y)\,\rme^{-\rmi y}\dd y~. 
\] 
(Remark that $P$ has zero mean over a period.) Replacing $P(y)$ 
with $\rme^{-\rmi\phi} P(y+\xi)$ if needed, we can assume that both 
$a$ and $b$ are real. If $P$ (hence also $W$) is small, we have $T 
\approx \pi$ and $\Psi \approx 0$, hence $k \approx 1$ and $\ell 
\approx 0$.  This determines uniquely the expansion of $P,k,\ell$ 
in powers of $a$ and $b$. Setting 
\begin{equation}\label{e:Wabdef} 
  W_{a,b}(x) \,=\, \rme^{\rmi\ell_{a,b}x}P_{a,b}(k_{a,b}x)~, \quad 
  x \in \R~, 
\end{equation} 
we obtain after straightforward calculations: 
\begin{eqnarray}\nonumber 
  \ell_{a,b}&=&{1 \over 4}(b^2-a^2) 
     + \rmO(a^4 + b^4)~,\\[1ex] \label{e:ellkPdef} 
  k_{a,b}&=&1 - {3 \over 4}(a^2+b^2) + \rmO(a^4 + b^4)~,\\[1ex]\nonumber 
  P_{a,b}(y)&=& a\,\rme^{-\rmi y} + b\,\rme^{\rmi y} 
    -\, {a^2 b \over 8}\,\rme^{-3\rmi y} -{a b^2 \over 8}\, 
     \rme^{3\rmi y} 
    + \rmO(|ab|(|a|^3 + |b|^3))~, 
\end{eqnarray} 
as $(a,b)\to(0,0)$. Notice also that the invariants $J,E$ have the 
following expressions: 
\begin{eqnarray*} 
  J &=& b^2-a^2 + \frac12(a^4-b^4) + \rmO(a^6 + b^6)~,\\ 
  E &=& a^2+b^2 -3a^2b^2 - \frac34(a^4+b^4) + \rmO(a^6 + b^6)~. 
\end{eqnarray*} 
With this parametrization, replacing $a$ with $-a$ or $b$ with $-b$ 
gives the same function $P$ up to a translation and a 
phase factor: 
\[ 
  P_{-a,b}(y) \,=\, -\rmi P_{a,b}(y+\pi/2)~, \quad 
  P_{-a,-b}(y) \,=\, - P_{a,b}(y) \,=\, P_{a,b}(y+\pi)~, \quad y \in \R~. 
\] 
It follows that $J,E$, hence also $k,\ell$, are even functions of 
$a$ and $b$. Similarly, $P_{b,a}(y) = \overline{P_{a,b}(y)}$. This 
conjugation leaves $E$ unchanged but reverses the sign of $J$, 
hence $k_{a,b}=k_{b,a}$ and $\ell_{a,b}=-\ell_{b,a}$. Therefore, 
using the symmetries of \reff{e:sgl}, we can restrict ourselves to 
the parameter region $\{b \ge a \ge 0\}$ without loss of 
generality. 
 
Two particular cases will play a special role in what follows. 
\begin{enumerate} 
\item {\bf (Cnoidal waves)} If $a = b$, then $\ell_{a,a} = 0$ and  
we obtain a family of real-valued periodic solutions $W_{a,a}(x) =  
P_{a,a}(k_{a,a}x)$, where 
\[ 
  k_{a,a} \,=\, 1 - {3 \over 2}\,a^2  + \rmO(a^4)~, \quad 
  P_{a,a}(y) \,=\, 2a \cos y   - {a^3 \over 4}\,\cos(3y) + \rmO(|a|^5)~. 
\] 
Observe that $J = 0$ in that case. 
 
\item {\bf (Plane waves)} If $a = 0$, then $P_{0,b}(y) = b\, 
\rme^{\rmi y}$, hence $W_{0,b}$ has constant modulus. It  
follows that 
\[ 
  W_{0,b}(x) \,=\, b\,\rme^{\rmi\sqrt{1-b^2}x}~, \quad \hbox{and} 
  \quad k_{0,b}+\ell_{0,b} \,=\, \sqrt{1-b^2}~. 
\] 
In addition, one also finds $k_{0,b}=\sqrt{1-3b^2/2}$. It is 
advantageous here to use the representation \reff{e:Qrep}, namely 
$W_{0,b}(y) = \rme^{\rmi p_{0,b}^+ y}Q_{0,b}^+(y)$ with $p_{0,b}^+ 
= \ell_{0,b}+k_{0,b} = (1-b^2)^{1/2}$ and $Q_{0,b}^+ \equiv b$. 
\\[2mm] 
Similarly, if $b=0$ we have $P_{a,0}(y)=a\rme^{-\rmi y}$ 
and thus 
\[ 
  W_{a,0}(x) \,=\, a\,\rme^{-\rmi\sqrt{1-a^2}x}~,\quad 
  \hbox{and}\quad k_{a,0}-\ell_{a,0}=\sqrt{1-a^2}~, \quad 
  k_{a,0} \,=\, \sqrt{1-3a^2/2}~. 
\] 
Alternatively, $W_{a,0}(y) = \rme^{\rmi p_{a,0}^- y}Q_{a,0}^-(y)$ 
with $p_{a,0}^- = \ell_{a,0}- k_{a,0} = -(1-a^2)^{1/2}$ and 
$Q_{a,0}^- \equiv a$. 
\end{enumerate}

\section{Orbital stability} 
\label{s:os} 
 
In this section we prove the orbital stability result in 
Theorem~\ref{th:orbit}. Since we restrict ourselves to periodic waves 
with small amplitude, we shall use the local parametrization 
\reff{e:Wabdef}, \reff{e:ellkPdef} of the small solutions of 
\reff{e:sgl}. Given $(a,b) \in \R^2$ with $\|(a,b)\|$ sufficiently 
small, we consider the periodic wave $U_{a,b}(x,t) = \rme^{-\rmi t} 
W_{a,b}(x)$, where 
\[ 
  W_{a,b}(x) \,=\, \rme^{\rmi \ell_{a,b}x} P_{a,b}(k_{a,b} x) \,=\, 
  \rme^{\rmi p_{a,b}x} Q_{a,b}(2k_{a,b}x)~, \quad x \in \R~. 
\] 
Here $\ell_{a,b}, k_{a,b}, P_{a,b}$ are defined in \reff{e:ellkPdef}, 
and the last expression in the right-hand side corresponds to the 
first choice in \reff{e:Qrep}, namely 
\begin{equation}\label{e:qQdef} 
  p_{a,b} \,=\, \ell_{a,b} + k_{a,b}~, \quad 
  Q_{a,b}(z) \,=\, \rme^{-\rmi z/2} P_{a,b}(z/2)~. 
\end{equation} 
From the properties of $P_{a,b}$ we deduce 
\begin{equation}\label{e:symQ} 
Q_{-a,b}(z)=Q_{a,b}(z+\pi),\quad Q_{-a,-b}(z)=-Q_{a,b}(z),\quad 
Q_{b,a}(z)=\rme^{-\rmi z} \overline{Q_{a,b}(z)}~, 
\end{equation} 
and that the real and imaginary parts of $Q_{a,b}$ are even and 
odd functions of $z$, respectively. 
 
\begin{Remark}\label{r:Qpm} 
Without loss of generality, we shall assume henceforth that $b \ge 
a\ge 0$. Note that the second choice in \reff{e:Qrep} would be 
preferable when $a^2\geq b^2$. 
\end{Remark} 
 
\subsection{Main result and strategy of proof} 
 
To study the stability of $U_{a,b}(x,t)$ we consider solutions of 
\reff{e:nls} of the form 
\begin{equation}\label{e:Qgen} 
  U(x,t) \,=\, \rme^{\rmi (p_{a,b}x -t)} Q(2k_{a,b}x,t)~, 
\end{equation} 
where $Q(z,t)$ is a $2\pi$-periodic function of $z$ which satisfies 
the evolution equation 
\begin{equation}\label{e:Qeq} 
  \rmi Q_t + 4\rmi p_{a,b}k_{a,b} Q_z + 4k_{a,b}^2 Q_{zz} + 
  (1-p_{a,b}^2)Q - |Q|^2Q \,=\, 0~. 
\end{equation} 
By construction, $Q_{a,b}(z)$ is now a stationary solution of 
\reff{e:Qeq} and our goal is to show that this equilibrium is 
stable with respect to $2\pi$-periodic perturbations. We thus 
introduce the function space 
\[ 
  X \,=\, H^1_\per([0,2\pi],\C) \,=\, \Bigl\{u \in H^1_\loc(\R,\C) 
  \,\Big|\,u(z{+}2\pi) = u(z)\,, ~\forall\, z \in \R\Bigr\}~, 
\] 
which is viewed as a {\it real} Hilbert space equipped with the 
scalar product 
\[ 
  (u,v)_X \,=\, \Re \int_0^{2\pi} (u(z) \overline v(z) + 
  u_z(z) \overline v_z(z)) \dd z~, \quad u,v \in X~. 
\] 
As usual, the dual space $X^*$ will be identified with 
$H^{-1}_\per([0,2\pi],\C)$ through the pairing 
\[ 
  \langle u,v \rangle \,=\, \Re \int_0^{2\pi} u(z) \overline v(z) \dd z~, 
  \quad u \in X^*~, \quad v \in X~. 
\] 
 
It is well-known that the Cauchy problem for \reff{e:Qeq} is globally 
well-posed in the space $X$. Moreover, the evolution defined by 
\reff{e:Qeq} on $X$ is invariant under a two-parameter group of 
isometries. The symmetry group is the two-dimensional torus 
$G = \T^2 = (\R / 2\pi\Z)^2$, a compact abelian Lie group which 
acts on $X$ through the unitary representation $\RR$ 
defined by 
\[ 
  (\RR_{(\phi,\xi)} u)(z) \,=\, \rme^{-\rmi \phi} u(z+\xi)~, 
  \quad u \in X~, \quad (\phi,\xi) \in G~. 
\] 
Due to these symmetries, it is useful to introduce the semi-distance 
$\rho$ on $X$ defined by 
\begin{equation}\label{e:rhodef} 
  \rho(u,v) \,=\, \inf_{(\phi,\xi) \in G} \|u - \RR_{(\phi,\xi)}v\|_X~, 
  \quad u,v \in X~. 
\end{equation} 
In words, $\rho(u,v)$ is small if $u$ is close to $v$ (in the 
topology of $X$) up to a translation and a phase rotation. Our 
stability result in Theorem~\ref{th:orbit} can now be formulated 
as follows: 
 
\begin{Proposition}\label{th:Qprop} 
There exist $C_0 > 0$, $\epsilon_0 > 0$, and $\delta_0 > 0$ such 
that, for all $(a,b) \in \R^2$ with $\|(a,b)\| \le \delta_0$, the 
following holds. If $Q_0 \in X$ satisfies $\rho(Q_0,Q_{a,b}) \le 
\epsilon$ for some $\epsilon \le \epsilon_0$, then the solution 
$Q(z,t)$ of \reff{e:Qeq} with initial data $Q_0$ satisfies 
$\rho(Q(\cdot,t), Q_{a,b}) \le C_0 \epsilon$ for all $t \in \R$. 
\end{Proposition} 
 
For each fixed value of $(a,b)$, the stability of the periodic (or 
plane) wave $Q_{a,b}$ can be proved using the abstract results of 
Grillakis, Shatah and Strauss \cite{GSS1,GSS2}. However, this 
approach would not give a stability theorem that holds uniformly 
in a neighborhood of the origin, as it is the case in 
Proposition~\ref{th:Qprop}. A difficulty in proving such a uniform 
result is that we have to deal simultaneously with three sorts of 
solutions: the zero solution ($a = b = 0$), plane waves ($ab = 0$) 
and periodic waves ($ab \neq 0$). These equilibria are genuinely 
different from the point of view of orbital stability theory, 
because their orbits under the action of the symmetry group $G$ 
have different dimensions ($0$, $1$, and $2$, respectively). In 
what follows, we shall concentrate on the periodic waves, and at 
the end we shall indicate how the other cases can be incorporated 
to obtain a uniform result. Whenever possible, we shall adopt 
similar notations as in \cite{GSS2} to facilitate comparison. 
 
Due to its symmetries, equation \reff{e:Qeq} has the same 
conserved quantities as the original NLS equation, namely 
\begin{eqnarray*} 
  && N(Q) \,=\, \frac12 \int_0^{2\pi} |Q(z)|^2 \dd z~, \\[2mm] 
  && M(Q) \,=\, \frac{\rmi}{2} \int_0^{2\pi} \overline Q(z) 
     Q_z(z) \dd z~, \\[2mm] 
  && \EE(Q) \,=\, \int_0^{2\pi} \Bigl( 2k_{a,b}^2 |Q_z(z)|^2 
     + \frac14|Q(z)|^4\Bigr)\dd z~. 
\end{eqnarray*} 
The charge $N$, the momentum $M$ and the energy $\EE$ are smooth 
real-valued functions on $X$. Their first order derivatives are 
therefore smooth maps from $X$ into $X^*$: 
\[ 
  N'(Q) \,=\, Q~, \quad M'(Q) \,=\, \rmi Q_z~, \quad 
  \EE'(Q) \,=\, -4k_{a,b}^2 Q_{zz} + |Q|^2 Q~. 
\] 
Similarly, the second order derivatives are smooth maps from 
$X$ into $\mathcal{L}(X,X^*)$, the space of all bounded linear 
operators from $X$ into $X^*$: 
\[ 
  N''(Q) = \1~, \quad M''(Q) \,=\, \rmi \partial_z~, \quad 
  \EE''(Q) \,=\, -4k_{a,b}^2\partial_{zz} + |Q|^2 + 2 Q\otimes Q~, 
\] 
where 
\[ 
\langle (Q\otimes Q)u,v\rangle = \int_0^{2\pi} 
\Re(Q\overline{u})\Re(Q\overline{v})\dd z, \quad \forall\ u,v \in 
X. 
\] 
 
From now on, we fix $(a,b) \in \R^2$ with $\|(a,b)\|$ sufficiently 
small. As is explained above, we assume for the moment that $ab 
\neq 0$, in which case the function $Q_{a,b} \in X$ defined by 
\reff{e:qQdef} is a stationary solution of \reff{e:Qeq} 
corresponding to a periodic wave of the original NLS equation, 
i.e. $|Q_{a,b}|$ is not constant. By construction, $Q_{a,b}$ is a 
critical point of the modified energy 
\begin{equation}\label{e:Eabdef} 
  \EE_{a,b}(Q) \,=\, \EE(Q) - (1-p_{a,b}^2)N(Q) - 4p_{a,b}k_{a,b}M(Q)~, 
\end{equation} 
namely $\EE_{a,b}'(Q_{a,b}) = 0$. The orbital stability argument 
is based on two essential ingredients: 
 
\medskip 
\noindent{\bf Claim 1:} The equilibrium $Q_{a,b}$ is a {\em local minimum} 
of the function $\EE_{a,b}$ restricted to the codimension two 
submanifold 
\begin{equation}\label{e:Sabdef} 
  \Sigma_{a,b}(Q) \,=\, \Big\{Q \in X\,\Big|\, N(Q) = N(Q_{a,b})\,,~ 
  M(Q) = M(Q_{a,b})\Big\}~. 
\end{equation} 
Note that this manifold contains the entire orbit of $Q_{a,b}$ 
under the action of $G$. 
 
\medskip\noindent{\bf Claim 2:} 
The equilibrium $Q_{a,b}$ is a member of a {\em two-parameter family} 
of travelling and rotating waves of the form 
\begin{equation}\label{e:travrot} 
  Q(z,t) \,=\, \rme^{-\rmi \omega t} \,Q_{a,b}^{\omega,c}(z+ct)~, 
  \quad z \in \R~, \quad t \in \R~, 
\end{equation} 
where $(\omega,c)$ lies in a neighborhood of the origin in $\R^2$ 
(the Lie algebra of $G$) and $Q_{a,b}^{\omega,c} \in X$ is a  
smooth function of $(\omega,c)$ with $Q_{a,b}^{0,0} = Q_{a,b}$.  
Moreover the map $(\omega,c) \mapsto (N(Q_{a,b}^{\omega,c}), 
M(Q_{a,b}^{\omega,c}))$ is a \emph{local diffeomorphism} near  
$(\omega,c) = (0,0)$. 
 
\subsection{Proof of Claim 2} 
 
The second claim is easily justified using the continuous 
symmetries of the NLS equation. Indeed, let $(a',b') \in \R^2$ be 
close to $(a,b)$. Then 
\[ 
  U(x,t) \,=\, \rme^{\rmi (p_{a'\!,b'}x -t)} 
  \,Q_{a'\!,b'}(2k_{a'\!,b'}x)~, \quad x \in \R~, \quad t \in \R~, 
\] 
is a solution of the NLS equation, but it is not of the form 
\reff{e:Qgen} because $p_{a'\!,b'} \neq p_{a,b}$ and $k_{a'\!,b'} 
\neq k_{a,b}$ in general. However we can transform $U(x,t)$ into a 
solution of \reff{e:nls} of the form \reff{e:Qgen}, 
\reff{e:travrot} by applying successively a dilation of factor 
$\lambda$ and a Galilean transformation of speed $v$, where 
\begin{equation}\label{e:lamvdef} 
  \lambda \,=\, \lambda_{a,b}^{a'\!,b'} \,=\, \frac{k_{a,b}}{k_{a'\!,b'}}~, 
  \quad v = v_{a,b}^{a'\!,b'} \,=\, 2(\lambda_{a,b}^{a'\!,b'}p_{a'\!,b'} - 
  p_{a,b})~. 
\end{equation} 
After some elementary algebra, we obtain $Q_{a,b}^{\omega,c}(z) = 
\lambda_{a,b}^{a'\!,b'} Q_{a'\!,b'}(z)$ with 
\begin{equation}\label{e:omcdef} 
  \omega \,=\, (\lambda_{a,b}^{a'\!,b'})^2(1 - p_{a'\!,b'}^2) - 
  (1 - p_{a,b}^2)~, \quad c \,=\, 4 (\lambda_{a,b}^{a'\!,b'})^2 
  k_{a'\!,b'}p_{a'\!,b'} - 4 k_{a,b}p_{a,b}~. 
\end{equation} 
Using the expansions \reff{e:ellkPdef}, it is straightforward to 
verify that 
\[ 
  \MM_{a,b} \,\eqdef\, \pmatrix{ 
    \frac{\partial \omega}{\partial a'} & 
    \frac{\partial c}{\partial a'} \vspace{1mm}\cr 
    \frac{\partial \omega}{\partial b'} & 
    \frac{\partial c}{\partial b'}} 
  \bigg|_{(a',b') = (a,b)} \,=\, 
  \pmatrix{4a & -2a \cr 2b & 2b}(\1 + \rmO(a^2{+}b^2))~. 
\] 
Since we assumed that $ab \neq 0$, the matrix $\MM_{a,b}$ is 
invertible for $a,b$ sufficiently small, hence the mapping 
$(a',b') \mapsto (\omega,c)$ defined by \reff{e:omcdef} 
is a diffeomorphism from a neighborhood of $(a,b)$ onto a 
neighborhood of $(0,0)$. This proves the existence of the travelling 
and rotating wave \reff{e:travrot} for $(\omega,c) \in \R^2$ 
sufficiently small. Remark that the profile $Q_{a,b}^{\omega,c}$ is a 
critical point of the functional 
\[ 
  \EE_{a,b}^{\omega,c}(Q) \,=\, \EE_{a,b}(Q) - \omega N(Q) - c M(Q)~, 
  \quad Q \in X~. 
\] 
Following \cite{GSS2}, we define $d_{a,b}(\omega,c) = 
\EE_{a,b}^{\omega,c}(Q_{a,b}^{\omega,c})$. The properties of the 
function $d_{a,b}$ will play an important role in the orbital 
stability argument. 
 
\begin{Lemma}\label{th:dlem} 
The Hessian matrix of the function $d_{a,b}$ satisfies: 
\[ 
  \HH_{a,b} \,\eqdef\, \pmatrix{ 
    \frac{\partial^2 d_{a,b}}{\partial \omega^2} & 
    \frac{\partial^2 d_{a,b}}{\partial \omega \,\partial c} \vspace{2mm}\cr 
    \frac{\partial^2 d_{a,b}}{\partial c \,\partial \omega} & 
    \frac{\partial^2 d_{a,b}}{\partial c^2}} 
  \Bigg|_{(\omega,c) = (0,0)} \,=\, \frac{\pi}3 
  \pmatrix{-2 & -1 \cr -1 & 1}(\1 + \rmO(a^2{+}b^2))~. 
\] 
\end{Lemma} 
 
\noindent{\bf Proof.} Since $Q_{a,b}^{\omega,c}$ is a critical point 
of $\EE_{a,b}^{\omega,c}$, we have 
\begin{equation}\label{e:dNM} 
  \frac{\partial}{\partial \omega}d_{a,b}(\omega,c) \,=\, 
    -N(Q_{a,b}^{\omega,c})~, \quad 
  \frac{\partial}{\partial c}d_{a,b}(\omega,c) \,=\, 
    -M(Q_{a,b}^{\omega,c})~. 
\end{equation} 
To compute the second-order derivatives, we parametrize $(\omega,c)$ 
by $(a',b')$ as above. Using \reff{e:omcdef} we find $\HH_{a,b} = 
-(\MM_{a,b})^{-1} \KK_{a,b}$, where 
\[ 
  \KK_{a,b} \,=\, \pmatrix{ 
    \frac{\partial}{\partial a'} N(Q_{a,b}^{\omega,c}) & 
    \frac{\partial}{\partial a'} M(Q_{a,b}^{\omega,c}) \vspace{2mm}\cr 
    \frac{\partial}{\partial b'} N(Q_{a,b}^{\omega,c}) & 
    \frac{\partial}{\partial b'} M(Q_{a,b}^{\omega,c})} 
  \Bigg|_{(a',b') = (a,b)}~. 
\] 
As $Q_{a,b}^{\omega,c} = \lambda_{a,b}^{a'\!,b'} Q_{a'\!,b'}(z)$, 
we have 
\[ 
  N(Q_{a,b}^{\omega,c}) \,=\, (\lambda_{a,b}^{a'\!,b'})^2 N(Q_{a'\!,b'})~, 
  \quad M(Q_{a,b}^{\omega,c}) \,=\, (\lambda_{a,b}^{a'\!,b'})^2 
  M(Q_{a'\!,b'})~. 
\] 
On the other hand, using the expansion 
\begin{equation}\label{e:Qexp} 
  Q_{a,b}(z) \,=\, a\,\rme^{-\rmi z} + b -{a^2 b \over 8} 
  \rme^{-2\rmi z} -{a b^2 \over 8}\rme^{\rmi z} 
  + \rmO(|ab|(|a|^3 + |b|^3))~, 
\end{equation} 
which follows from \reff{e:ellkPdef}, \reff{e:qQdef}, we easily find 
\[ 
  N(Q_{a,b}) \,=\, \pi(a^2 + b^2) + \rmO(a^2b^2(a^2+b^2))~, 
  \quad M(Q_{a,b}) \,=\, \pi a^2 + \rmO(a^2b^2(a^2+b^2))~. 
\] 
Combining these results, we obtain 
\[ 
  (\MM_{a,b})^{-1} \,=\, \frac{1}{6ab}\pmatrix{b & a \cr -b & 2a} 
 (\1 + \rmO(a^2{+}b^2))~, \quad \KK_{a,b} \,=\, 2\pi 
 \pmatrix{a & a \cr b & 0}(\1 + \rmO(a^2{+}b^2))~, 
\] 
and the conclusion follows. \QED 
 
\medskip 
Lemma~\ref{th:dlem} implies that the Hessian matrix $\HH_{a,b}$ is 
nondegenerate for $\|(a,b)\|$ sufficiently small (in fact, $\HH_{a,b}$ 
has one positive and one negative eigenvalue.) It follows that the map 
$(\omega,c) \mapsto (N(Q_{a,b}^{\omega,c}),M(Q_{a,b}^{\omega,c}))$ is 
a local diffeomorphism near $(\omega,c) = (0,0)$, because by 
\reff{e:dNM} the Jacobian matrix of this map at the origin is just 
$-\HH_{a,b}$. Thus Claim~2 above is completely justified. 
 
\begin{Remark}\label{r:gss} 
At this point we could apply the general result of \cite{GSS2}, 
but as already mentioned this would not give the uniform result in 
Theorem~\ref{th:orbit}. According to the {\it Stability Theorem} 
in \cite{GSS2}, in order to establish the orbital stability of a 
single wave $Q_{a,b}$ with $ab \neq 0$ it suffices to show that 
the linear operator 
\begin{equation}\label{e:Habdef} 
  H_{a,b} \,=\, \EE_{a,b}''(Q_{a,b}) \,=\, -4k_{a,b}^2 \partial_{zz} 
  -4\rmi p_{a,b}k_{a,b}\partial_z - (1{-}p_{a,b}^2) 
  + |Q_{a,b}|^2 + 2 Q_{a,b} \otimes Q_{a,b}~, 
\end{equation} 
has precisely one simple negative eigenvalue, a two-dimensional 
kernel spanned by 
\begin{equation}\label{e:kernel} 
  \frac{\partial}{\partial \phi} \RR_{(\phi,\xi)} 
  Q_{a,b}\Big|_{(\phi,\xi)=(0,0)} \,=\, - \rmi Q_{a,b} ~, \quad 
  \frac{\partial}{\partial \xi} \RR_{(\phi,\xi)} 
  Q_{a,b}\Big|_{(\phi,\xi)=(0,0)}\,=\, \partial_z Q_{a,b} ~, 
\end{equation} 
and that the rest of its spectrum is strictly positive. Observe 
that $H_{a,b}$ is \emph{self-adjoint} in the real Hilbert space 
$L^2_\per([0,2\pi],\C)$ equipped with the scalar product $\langle 
\cdot,\cdot \rangle$. Clearly, the vectors \reff{e:kernel} always 
belong to the kernel of $H_{a,b}$. In fact, for small $(a,b)$ we 
can determine the spectrum of $H_{a,b}$ by a perturbation argument 
similar to the one used for the spectral analysis of the operators 
$\mathcal A_{a,b,\gamma}$ in Section~\ref{s:ss}. We find that 
$H_{a,b}$ has exactly four eigenvalues in a neighborhood of the 
origin, the rest of the spectrum being positive and bounded 
away from zero.  Among these four eigenvalues, two are always 
zero, and the other two have negative product $-12a^2b^2(1 + 
\rmO(a^2{+}b^2))$. This implies that $H_{a,b}$ has the required 
properties, so that the wave profile $Q_{a,b}$ is orbitally stable 
if $ab \neq 0$. This information on the spectrum of $H_{a,b}$ will 
not be used in the remainder of this section. However, since it 
provides the starting point for the stability analysis of large 
waves in \cite{GH2}, we give a brief proof in the Appendix. 
\end{Remark} 
 
\subsection{Proof of Claim 1} 
 
We now turn back to Claim~1 and study the behavior of the energy 
$\EE_{a,b}$ on the manifold $\Sigma_{a,b}$ defined by 
\reff{e:Sabdef}. In the arguments below, we assume $b \ge a > 0$, 
so exclude for the moment the plane wave corresponding to 
$a = 0$. Let $\TT_{a,b}$ be the tangent space to $\Sigma_{a,b}$ 
at the point $Q_{a,b}$: 
\[ 
  \TT_{a,b} \,=\,\Bigl\{Q \in X \,\Big|\, \langle N'(Q_{a,b}),Q\rangle 
  =  \langle M'(Q_{a,b}),Q\rangle = 0\Bigr\}~. 
\] 
Then $X = \TT_{a,b} \oplus \NN_{a,b}$, where $\NN_{a,b}$ (the 
``normal'' space) is the two-dimensional subspace of $X$ spanned by 
$N'(Q_{a,b}) = Q_{a,b}$ and $M'(Q_{a,b}) = \rmi \partial_z Q_{a,b}$. 
When $(a,b)$ is small, a more convenient basis of $\NN_{a,b}$ is 
$\{\xi_{a,b},\eta_{a,b}\}$, where 
\begin{equation}\label{e:xietadef} 
  \xi_{a,b} \,=\, \frac{\rmi}{a}\partial_z Q_{a,b} \,=\, \rme^{-\rmi z} 
  + \rmO(|ab|+b^2)~, \quad 
  \eta_{a,b} \,=\, \frac{1}{b}(Q_{a,b} - \rmi \partial_z Q_{a,b}) 
  \,=\, 1 + \rmO(a^2+|ab|)~. 
\end{equation} 
The tangent space is further decomposed as $\TT_{a,b} = Y_{a,b} \oplus 
Z_{a,b}$, where 
\[ 
  Y_{a,b} \,=\,\Bigl\{Q \in \TT_{a,b} \,\Big|\, \langle \rmi Q_{a,b}, 
  Q\rangle =  \langle \partial_z Q_{a,b},Q\rangle = 0\Bigr\}~, 
\] 
and $Z_{a,b}$ is the two-dimensional space spanned by $\rmi Q_{a,b}$ 
and $\partial_z Q_{a,b}$. In view of \reff{e:kernel}, $Z_{a,b}$ is 
just the tangent space to the orbit of $Q_{a,b}$ under the action 
of $G$. Again, a convenient basis of $Z_{a,b}$ is $\{\rmi \xi_{a,b}, 
\rmi \eta_{a,b}\}$. 
 
As in \cite{GSS2}, we introduce an appropriate coordinate system in 
a neighborhood of the orbit of $Q_{a,b}$ under the action of $G$: 
 
\begin{Lemma}\label{th:coord} 
Assume that $\|(a,b)\|$ is sufficiently small and $b \ge a > 0$. 
There exist $\kappa > 0$, $C_1 > 0$, and $C_2 > 0$ such that any 
$Q \in X$ with $\rho(Q,Q_{a,b}) \le \kappa a$ can be represented 
as 
\begin{equation}\label{e:Qdecomp} 
  Q \,=\, \RR_{(\phi,\xi)}(Q_{a,b} + \nu + y)~, 
\end{equation} 
where $(\phi,\xi) \in G$, $\nu \in \NN_{a,b}$, $y \in Y_{a,b}$, and 
$\|\nu\|_X + \|y\|_X \le C_1 \rho(Q,Q_{a,b})$. Moreover, if $Q \in 
\Sigma_{a,b}$, then $\|\nu\|_X \le (C_2/a)\|y\|_X^2$. 
\end{Lemma} 
 
\begin{Remark} 
Here and in the sequel, all constants $C_1, C_2, \dots$ are 
independent of $(a,b)$ provided $\|(a,b)\|$ is sufficiently small. 
\end{Remark} 
 
\noindent{\bf Proof.} It is clearly sufficient to prove 
the result for all $Q \in X$ with $\|Q-Q_{a,b}\|_X \le \kappa a$, 
where $\kappa > 0$ is a (small) constant that will be fixed below. 
Since $X = \NN_{a,b} \oplus Y_{a,b} \oplus Z_{a,b}$, any such $Q$ 
can be decomposed as 
\[ 
  Q \,=\, Q_{a,b} -x_1 \rmi Q_{a,b} + x_2 \partial_z Q_{a,b} 
  + \nu_1 + y_1~, 
\] 
where $\nu_1 \in \NN_{a,b}$, $y_1 \in Y_{a,b}$, and $(x_1,x_2) \in 
\R^2$ is the solution of the linear system 
\begin{equation}\label{e:linsys} 
  \pmatrix{ 
  \langle \rmi Q_{a,b},\rmi Q_{a,b} \rangle & 
  -\langle \rmi Q_{a,b},\partial_z Q_{a,b} \rangle \cr 
  -\langle \partial_z Q_{a,b},\rmi Q_{a,b} \rangle & 
  \langle \partial_z Q_{a,b},\partial_z Q_{a,b} \rangle} 
  \pmatrix{x_1 \cr x_2} \,=\, 
  \pmatrix{-\langle \rmi Q_{a,b}, Q - Q_{a,b} \rangle \cr 
           \langle \partial_z Q_{a,b}, Q - Q_{a,b} \rangle}~. 
\end{equation} 
The matrix $\hat\MM_{a,b}$ in the left-hand side of 
\reff{e:linsys} is invertible, and using the expansions 
\reff{e:Qexp} we find 
\[ 
  (\hat\MM_{a,b})^{-1} \,=\, \frac{1}{2\pi}\pmatrix{b^{-2} & -b^{-2} 
  \cr -b^{-2}  & a^{-2} + b^{-2}}(\1 + \rmO(a^2{+}b^2))~. 
\] 
Since $b\geq a$, it follows that $|x_1| + |x_2| \le (C/a)\|Q - 
Q_{a,b}\|_X \le C\kappa$ for some $C > 0$ (independent of $a,b$). 
Now, since 
\[ 
  \RR_{(\phi,\xi)} Q_{a,b} \,=\, Q_{a,b} -\phi \rmi Q_{a,b} + \xi 
  \partial_z Q_{a,b} + \rmO(\phi^2 + \xi^2)~, 
\] 
the Implicit Function Theorem implies that, if $(x_1,x_2)$ is 
sufficiently small, there exists a unique pair $(\phi,\xi) \in 
\R^2$ with $(\phi,\xi) = (x_1,x_2) + \rmO(x_1^2+x_2^2)$ such that 
$\RR_{(\phi,\xi)}^{-1} Q - Q_{a,b} \in \NN_{a,b} \oplus Y_{a,b}$ 
(see Lemma~4.2 in \cite{GSS2} for a similar argument). Setting 
$\RR_{(\phi,\xi)}^{-1} Q - Q_{a,b} = \nu + y$, we obtain the 
desired decomposition (assuming that $\kappa > 0$ is small enough 
so that we can apply the Implicit Function Theorem). This choice 
of $(\phi,\xi)$ does not minimize the distance in $X$ between $Q$ 
and $\RR_{(\phi,\xi)} Q_{a,b}$, because the subspaces $\NN_{a,b}$, 
$Z_{a,b}$, and $Y_{a,b}$ are not mutually orthogonal for the 
scalar product of $X$. However, since the \emph{minimum gap} 
between these spaces is strictly positive (uniformly in $a,b$), we 
still have $\|\nu\|_X + \|y\|_X \le C \|\nu + y\|_X \le C_1 
\rho(Q,Q_{a,b})$. (We refer to \cite{K1} for the definition and 
the properties of the minimum gap between closed subspaces of a 
Banach space.) 
 
Now, we assume in addition that $Q \in \Sigma_{a,b}$, i.e. 
$N(Q) = N(Q_{a,b})$ and $M(Q) = M(Q_{a,b})$. In view of 
\reff{e:Qdecomp}, we have 
\[ 
  N(Q) \,=\, N(Q_{a,b} + \nu + y) \,\equiv\, N(Q_{a,b}) + 
  \langle Q_{a,b} , \nu + y \rangle + N(\nu + y)~, 
\] 
and using the fact that $y \in Y_{a,b} \subset \TT_{a,b}$ we 
obtain $\langle Q_{a,b} , \nu \rangle + N(\nu + y) = 0$. 
A similar argument shows that $\langle \rmi \partial_z Q_{a,b} , 
\nu \rangle + M(\nu + y) = 0$. Thus $\nu = \nu_1 Q_{a,b} + \nu_2 
\rmi \partial_z Q_{a,b}$, where 
\[ 
  \pmatrix{ 
  \langle Q_{a,b},Q_{a,b} \rangle & 
  \langle Q_{a,b},\rmi \partial_z Q_{a,b} \rangle \cr 
  \langle \rmi \partial_z Q_{a,b},Q_{a,b} \rangle & 
  \langle \rmi \partial_z Q_{a,b},\rmi \partial_z Q_{a,b} \rangle} 
  \pmatrix{\nu_1 \cr \nu_2} \,=\, - 
  \pmatrix{N(\nu+y) \cr M(\nu+y)}~. 
\] 
Observe that the matrix of this system is exactly the same one 
as in \reff{e:linsys}. Thus, proceeding as above, we obtain $\|\nu\|_X 
\le (C/a)\|\nu + y\|_X^2$ for some $C > 0$ independent of $a,b$. 
Since we already know that $\|\nu\|_X \le C_1 \kappa a$, it 
follows that $\|\nu\|_X \le (C_2/a)\|y\|_X^2$ provided $\kappa 
> 0$ is sufficiently small. \QED 
 
\medskip 
To show that the energy $\EE_{a,b}$ has a local minimum on 
$\Sigma_{a,b}$ at $Q_{a,b}$, we consider the second variation 
of $\EE_{a,b}$ at $Q_{a,b}$, i.e. the linear operator $H_{a,b}$ 
defined in \reff{e:Habdef}. 
 
\begin{Lemma}\label{th:positive} 
If $\|(a,b)\|$ is sufficiently small and $ab \neq 0$, then 
\begin{equation}\label{e:positive} 
  \langle H_{a,b} \,y, y\rangle \,\ge\, 6 \|y\|_X^2~, \quad 
  \hbox{for all } y \in Y_{a,b}~. 
\end{equation} 
\end{Lemma} 
 
\noindent{\bf Proof.} We use a perturbation argument. When $(a,b) 
= (0,0)$, the operator $H_{a,b}$ reduces to a differential 
operator with constant coefficients: $H_0 = -4\partial_{zz} 
-4\rmi\partial_z$. This operator is self-adjoint in the real 
Hilbert space $L^2_\per([0,2\pi],\C)$ equipped with the scalar 
product $\langle \cdot,\cdot \rangle$, and its spectrum is 
$\sigma(H_0) = \{4n(n\pm1) \,|\, n \in \Z\}$. The null space of 
$H_0$ is spanned by the four vectors $\xi_0, \rmi \xi_0, \eta_0, 
\rmi \eta_0$, where $\xi_0 = \rme^{-\rmi z}$ and $\eta_0 = 1$ (see 
\reff{e:xietadef}). The other eigenvalues of $H_0$ are positive 
and greater or equal to $8$, hence the quadratic form $h_0 : X \to 
\R$ associated to $H_0$ satisfies 
\[ 
  h_0(Q) \,\eqdef\, \langle H_0 Q,Q \rangle \,\ge\, 8 \|Q\|_X^2~, 
  \quad \hbox{for all } Q \in Y_0~, 
\] 
where 
\[ 
  Y_0 \,=\, \Bigl\{Q \in X\,|\, \langle \xi_0,Q \rangle = 
  \langle \rmi \xi_0,Q \rangle = \langle \eta_0,Q \rangle = 
  \langle \rmi \eta_0,Q \rangle = 0 \Bigr\}~. 
\] 
We now consider the quadratic form $h_{a,b} : X \to \R$ defined by 
$h_{a,b}(Q) = \langle H_{a,b} Q,Q \rangle$. This form is uniformly 
bounded for $(a,b)$ in a neighborhood of zero, i.e. there exists 
$C_3 > 0$ such that $h_{a,b}(Q) \le C_3 \|Q\|_X^2$ for all $Q \in X$. 
Moreover, $h_{a,b}$ converges to $h_0$ as $(a,b) \to (0,0)$ in the 
following sense: 
\[ 
  \sup\Bigl\{|h_{a,b}(Q) - h_0(Q)| \,\Big|\, Q \in X\,,~ 
  \|Q\|_X = 1\Bigr\} \,=\, \rmO(a^2+b^2)~. 
\] 
In particular, we have $h_{a,b}(Q) \ge 7\|Q\|_X^2$ for all $Q \in 
Y_0$ if $\|(a,b)\|$ is sufficiently small. 
 
On the other hand, since $\|\xi_{a,b} - \xi_0\|_X + \|\eta_{a,b} 
- \eta_0\|_X = \rmO(a^2+b^2)$, it is straightforward to verify that 
the subspace 
\[ 
  Y_{a,b} \,=\, \Bigl\{Q \in X\,|\, \langle \xi_{a,b},Q \rangle = 
  \langle \rmi \xi_{a,b},Q \rangle = \langle \eta_{a,b},Q \rangle = 
  \langle \rmi \eta_{a,b},Q \rangle = 0 \Bigr\} 
\] 
converges to $Y_0$ as $(a,b) \to (0,0)$ in the following sense: 
\[ 
  \delta(Y_{a,b},Y_0) \,\eqdef\, \sup\Bigl\{\dist_X(Q,Y_0) 
  \,\Big|\, Q \in Y_{a,b}\,,~ \|Q\|_X = 1\Bigr\} 
  \,=\, \rmO(a^2+b^2)~. 
\] 
In particular, if $Q \in Y_{a,b}$ satisfies $\|Q\|_X = 1$, we can 
find $\tilde Q \in Y_0$ with $\|\tilde Q\|_X = 1$ and $\|Q - 
\tilde Q\|_X$ as small as we want, provided $(a,b)$ is close to 
zero. Since $h_{a,b}(\tilde Q) \ge 7$ and 
\[ 
  |h_{a,b}(Q) - h_{a,b}(\tilde Q)| \,\le\, \|h_{a,b}\| 
  (\|Q\|_X + \|\tilde Q\|_X)\|Q - \tilde Q\|_X 
  \,\le\, 2 C_3 \|Q - \tilde Q\|_X~, 
\] 
we conclude that $h_{a,b}(Q) \ge 6$ if $\|(a,b)\|$ is sufficiently 
small. This proves \reff{e:positive}. \QED 
 
\medskip 
Using Lemmas~\ref{th:coord} and \ref{th:positive}, we are able to  
give a more precise version of Claim~1 above: 
 
\begin{Lemma}\label{th:coercive} 
There exists $\kappa_1 > 0$ such that, if $\|(a,b)\|$ is 
sufficiently small and $b \ge a > 0$, then for all $Q \in 
\Sigma_{a,b}$ satisfying $\rho(Q,Q_{a,b}) \le \kappa_1 a$ one has 
the inequality 
\begin{equation}\label{e:coercive} 
  \EE_{a,b}(Q) - \EE_{a,b}(Q_{a,b}) \,\ge\, \rho(Q,Q_{a,b})^2~. 
\end{equation} 
\end{Lemma} 
 
\noindent{\bf Proof.} If $Q \in \Sigma_{a,b}$ satisfies $\rho(Q,Q_{a,b}) 
\le \kappa_1 a$ for some $\kappa_1 > 0$ sufficiently small, we know from 
Lemma~\ref{th:coord} that $Q = \RR_{(\phi,\xi)} (Q_{a,b} + \nu + y)$, 
where $(\phi,\xi) \in G$, $\nu \in \NN_{a,b}$, $y \in Y_{a,b}$, 
$\|y\|_X \le C_1 \rho(Q,Q_{a,b})$, and $\|\nu\|_X \le (C_2/a) 
\|y\|_X^2$. In particular, $\|\nu\|_X \le \kappa_1 C_1 C_2 \|y\|_X$. 
Since the energy $\EE_{a,b}$ is invariant under the 
action of $G$, we have $\EE_{a,b}(Q) = \EE_{a,b}(Q_{a,b} + \nu + y)$. 
As $\EE_{a,b}'(Q_{a,b}) = 0$ and $\EE_{a,b}''(Q_{a,b}) = H_{a,b}$, 
we obtain using Taylor's formula: 
\[ 
  \EE_{a,b}(Q) - \EE_{a,b}(Q_{a,b}) \,=\, \frac12 
  \langle H_{a,b}(y+\nu),(y+\nu)\rangle + \rmO(\|y\|_X^3)~. 
\] 
But $\langle H_{a,b} \,y,y\rangle \ge 6 \|y\|_X^2$ by 
Lemma~\ref{th:positive}, hence 
\begin{eqnarray*} 
  \frac12 \langle H_{a,b}(y+\nu),(y+\nu)\rangle &\ge& 
  3\|y\|_X^2 - C_3 \|y\|_X \|\nu\|_X -\half C_3 \|\nu\|_X^2 \\ 
  &\ge& \|y\|_X^2(3 - \kappa_1 C_1 C_2 C_3 - \half 
    (\kappa_1 C_1 C_2)^2 C_3)~, 
\end{eqnarray*} 
where $C_3$ is the constant in the proof of Lemma~\ref{th:positive}. 
Thus, if $\kappa_1$ is sufficiently small, we obtain $\EE_{a,b}(Q) - 
\EE_{a,b}(Q_{a,b}) \ge 2\|y\|_X^2$. Under the same assumption, we 
also have $\rho(Q,Q_{a,b}) \le \|y+\nu\|_X \le \|y\|_X (1+\kappa_1 
C_1 C_2) \le \sqrt{2}\|y\|_X$, and \reff{e:coercive} follows. \QED 
 
\subsection{Proof of Proposition~\ref{th:Qprop}} 
 
The proof of Proposition~\ref{th:Qprop} consists of three steps in 
which we treat successively the three types of waves: the zero 
solution ($a=b=0$), the plane waves ($ab=0$), and the periodic 
waves ($ab\not=0$). In each case, the arguments rely upon energy 
estimates as the one given in Lemma~\ref{th:coercive} for the 
periodic waves ($ab\not=0$). In the case of the plane wave 
$Q_{0,b} \equiv b$ the stability is proved using the functional 
\[ 
  \EE_b(Q) \,=\, \EE(Q) - b^2 N(Q) \,=\, 
  \int_0^{2\pi} \Bigl((2{-}3b^2)|Q_z|^2 + \frac14 (|Q|^2-b^2)^2 
 \Bigr)\dd z - \frac{\pi b^4}{2}~, 
\] 
for which we have the analog of Lemma~\ref{th:coercive}: 
 
\begin{Lemma}\label{th:coercive2} 
There exists $\kappa_2 > 0$ such that, if $b > 0$ is sufficiently  
small, then for all $Q \in X$ satisfying $\rho(Q,Q_{0,b}) \le \kappa_2 
b$ and $N(Q) = N(Q_{0,b})$ one has the inequality 
\begin{equation}\label{e:coercive2} 
  \EE_{b}(Q) - \EE_{b}(Q_{0,b}) \,\ge\, \frac16 \,\rho(Q,Q_{0,b})^2~. 
\end{equation} 
\end{Lemma} 
 
\begin{Proof} Without loss of generality, we assume that  
$\|Q-Q_{0,b}\|_X \le \kappa_2 b$. Since 
\[ 
  |Q(z)-b| \leq  \|Q-Q_{0,b}\|_X \le \kappa_2 b < b~,\quad \mbox{for 
  all}\ z\in\R~, 
\] 
we can write $Q = (b+r)\rme^{\rmi \phi}$, where $r,\phi \in X$ 
are real functions satisfying 
\[ 
  |r(z)|\leq \kappa_2 b~,\quad |\rme^{\rmi\varphi(z)}-1|\leq 2\kappa_2~, 
  \quad \mbox{for all}\ z\in\R~. 
\] 
These inequalities imply $\|r\|_X + b\|\phi\|_X \le C\kappa_2 b$, for 
some $C > 0$ independent of $a,b$. Thus 
\begin{eqnarray}\nonumber 
  \EE_{b}(Q) - \EE_{b}(Q_{0,b}) &=& \int_0^{2\pi} 
  \Bigl\{(2{-}3b^2)\Bigl(r_z^2 + (b+r)^2\phi_z^2\Bigr) 
  + \frac14 (2br+r^2)^2 \Bigr\}\dd z \\ \label{e:Aest} 
  & \ge & \int_0^{2\pi} (r_z^2 + b^2 \phi_z^2)\dd z~. 
\end{eqnarray} 
On the other hand, since $N(Q) = N(Q_{0,b})$, we have 
$\int_0^{2\pi} (2br+r^2)\dd z = 0$, and Poincar\'e's inequality 
implies 
\[ 
  \int_0^{2\pi} r^2 \Bigl(1 + \frac{r}{2b}\Bigr)^2\dd z 
  \,\le\, \int_0^{2\pi} r_z^2 \Bigl(1 + \frac{r}{b}\Bigr)^2\dd z~, 
  \quad \hbox{hence} \quad \int_0^{2\pi} r^2\dd z \,\le\, 
  2 \int_0^{2\pi} r_z^2\dd z~. 
\] 
Finally, if $\bar \phi = (2\pi)^{-1}\int_0^{2\pi} \phi(z)\dd z$, 
we have 
\begin{eqnarray}\nonumber 
  \rho(Q_,Q_{0,b})^2 & \le & \|Q - b\,\rme^{\rmi \bar\phi}\|_X^2 
  \,\le\, 2 \|r \rme^{\rmi\varphi}\|_X^2 + 2 
  b^2 \|\rme^{\rmi\varphi} - \rme^{\rmi\bar\varphi}\|_X^2\\ 
  &\leq & 2 \|r\|_X^2 + 2 \int_0^{2\pi} r^2\varphi_z^2 \dd z 
  + 2 b^2 \|\varphi - \bar\varphi\|_X^2 \,\leq\, 
  6 \int_0^{2\pi} (r_z^2 + b^2 \phi_z^2)\dd z~, \label{e:Best} 
\end{eqnarray} 
again by Poincar\'e's inequality. Combining \reff{e:Aest} and 
\reff{e:Best}, we obtain \reff{e:coercive2}~. 
\end{Proof} 
 
We are now in position to prove Proposition~\ref{th:Qprop}. 
 
\medskip\noindent {\bf Proof of Proposition~\ref{th:Qprop}.} 
Throughout the proof, we assume that $\|(a,b)\|$ is sufficiently 
small and that $b \ge a \ge 0$. Given $Q_0 \in X$ with $\rho(Q_0, 
Q_{a,b}) \le \epsilon$, we consider the solution $Q(z,t)$ of 
\reff{e:Qeq} with initial data $Q_0$. Replacing $Q_0$ with 
$\RR_{(\phi,\xi)} Q_0$ if needed, we can assume that $\|Q_0 - 
Q_{a,b}\|_X \le \epsilon$. We distinguish three cases: 
 
\smallskip\noindent{\bf Case 1:} $a = b = 0$, i.e. $Q_{a,b} = 0$. 
In this case, if $\epsilon > 0$ is small enough, the solution 
$Q(\cdot,t)$ of \reff{e:Qeq} satisfies $\|Q(\cdot,t)\|_X \le 
2\epsilon$ for all $t \in \R$. This is obvious because the quantity 
\[ 
  \EE(Q) + 4 N(Q) \,=\, \int_0^{2\pi} \Bigl(2|Q_z|^2 + 
  \frac14 |Q|^4 + 2|Q|^2\Bigr)\dd z 
\] 
is invariant under the evolution of \reff{e:Qeq}, and satisfies 
$2\|Q\|_X^2 \le \EE(Q) + 4 N(Q) \le 4\|Q\|_X^2$ if $\|Q\|_X$ is 
small. 
 
\noindent{\it Remark:} As a consequence of this preliminary case, 
we assume from now on that $\epsilon \le \kappa_3 (a^2+b^2)^{1/2}$ 
for some small $\kappa_3 > 0$. Indeed, if $\epsilon \ge \kappa_3 
(a^2+b^2)^{1/2}$, we can use the trivial estimate 
\[ 
  \|Q(\cdot,t) - Q_{a,b}\|_X \,\le\, \|Q(\cdot,t)\|_X + 
  \|Q_{a,b}\|_X \,\le\, 2\|Q_0\|_X + \|Q_{a,b}\|_X 
  \le 2\epsilon + 3\|Q_{a,b}\|_X~, 
\] 
which gives the desired result since $\|Q_{a,b}\|_X \le 
C(a^2+b^2)^{1/2} \le (C/\kappa_3)\epsilon$. 
 
\smallskip\noindent{\bf Case 2:} 
$b > a = 0$, i.e. $Q_{a,b} = b$ is a plane wave. 
We consider initial data $Q_0 \in X$ such that $\|Q_0 - 
Q_{0,b}\|_X \le \epsilon \,\le\, \kappa_3 b$. If $N(Q_0)$ = 
$N(Q_{0,b})$, then $N(Q(\cdot,t)) = N(Q_{0,b})$ for all $t \in 
\R$, and Lemma~\ref{th:coercive2} implies 
\[ 
  \rho(Q(\cdot,t),Q_{0,b})^2 \,\le\, 6(\EE_{b}(Q(\cdot,t)) - 
  \EE_{b}(Q_{0,b})) \,=\, 6(\EE_{b}(Q_0) - \EE_{b}(Q_{0,b})) 
  \,\le\, C \epsilon^2~, 
\] 
provided that $C\epsilon^2 \le \kappa_2^2 b^2$, which is the 
case if $C\kappa_3^2 \le \kappa_2^2$. If $N(Q_0) \neq 
N(Q_{0,b})$, we define $\omega = \pi^{-1}(N(Q_0) - N(Q_{0,b}))$, 
so that $N(Q_0) = N(Q_{0,b}^\omega)$, where $Q_{0,b}^\omega = 
(b^2+\omega)^{1/2}$. So we are led to study the stability of 
the rotating wave $Q_{0,b}^\omega \,\rme^{-\rmi \omega t}$ of 
\reff{e:Qeq} with respect to perturbations preserving the 
charge $N$. This can be proved exactly as above, and we 
obtain $\rho(Q(\cdot,t),Q_{0,b}^\omega) \le C\epsilon$ for 
all $t \in \R$. Since $\|Q_{0,b}^\omega - Q_{0,b}\|_X \le C 
|\omega|/b \le C\epsilon$, we have the desired result. 
 
\noindent{\it Remark:} As a consequence, we can assume 
from now on that $\epsilon \le \kappa_4 a$ for some small 
$\kappa_4 > 0$. Indeed, if $\epsilon \ge \kappa_4 a$, we can 
use the easy estimate 
\begin{equation}\label{e:easy} 
  \rho(Q(\cdot,t),Q_{a,b}) \,\le\, \rho(Q(\cdot,t),Q_{0,b}) 
  + \|Q_{0,b} - Q_{a,b}\|_X~. 
\end{equation} 
Observe that $\|Q_{0,b} - Q_{a,b}\|_X \le Ca$ for some $C > 0$ 
independent of $a,b$. In particular 
\[ 
  \|Q_0 - Q_{0,b}\|_X \,\le\, \|Q_0 - Q_{a,b}\|_X + 
  \|Q_{a,b} - Q_{0,b}\|_X \,\le\, \epsilon + Ca \,\le\, 
  (1+C/\kappa_4)\epsilon~, 
\] 
hence for $\epsilon > 0$ small enough we have $\rho(Q(\cdot,t), 
Q_{0,b}) \le C' \epsilon$ for all $t \in \R$. It then follows 
from \reff{e:easy} that $\rho(Q(\cdot,t),Q_{a,b}) \le C'' \epsilon$ 
for all $t \in \R$, which is the desired result. 
 
\smallskip\noindent{\bf Case 3:} $b \ge a > 0$, i.e. $Q_{a,b}$ 
is a nontrivial periodic equilibrium of \reff{e:Qeq} corresponding 
to a periodic wave of \reff{e:nls}. Assume that $Q_0 \in X$ 
satisfies $\|Q_0-Q_{a,b}\|_X \le \epsilon \le \kappa_4 a$. If 
$Q_0 \in \Sigma_{a,b}$, then $Q(\cdot,t) \in \Sigma_{a,b}$ for 
all $t \in \R$ and Lemma~\ref{th:coercive} implies that 
\[ 
  \rho(Q(\cdot,t),Q_{a,b})^2 \,\le\, \EE_{a,b}(Q(\cdot,t)) - 
  \EE_{a,b}(Q_{a,b}) \,=\, \EE_{a,b}(Q_0) - \EE_{a,b}(Q_{a,b}) 
  \,\le\, C\epsilon^2~, 
\] 
provided $C\epsilon^2 \le \kappa_1^2 a^2$, which is the case if 
$C\kappa_4^2 \le \kappa_1^2$.  If $Q_0 \notin \Sigma_{a,b}$, then by 
Claim~2 above there exists $(\omega,c) \in \R^2$ with $|\omega| + |c| 
\le Cb\epsilon $ such that $N(Q_0) = N(Q_{a,b}^{\omega,c})$ and 
$M(Q_0) = M(Q_{a,b}^{\omega,c})$. So we are led to study the stability 
of the periodic wave $u(x,t) = \,\rme^{\rmi (p_{a,b}x - 
(1+\omega)t)}\,Q_{a,b}^{\omega,c}(2k_{a,b}x+ct)$ of \reff{e:nls} 
among solutions of the form $\,\rme^{\rmi (p_{a,b}x-t)} 
\,Q(2k_{a,b}x,t)$ for which the charge $N$ and the momentum $M$ have 
the same values as for the periodic wave. But if we apply a dilation 
of factor $\lambda$ and a Galilean transformation of speed $v$, where 
$\lambda, v$ are given by \reff{e:lamvdef}, the periodic wave becomes 
$u(x,t) = \,\rme^{\rmi (p_{a'\!,b'}x -t)}\,Q_{a'\!,b'} 
(2k_{a'\!,b'}x)$ for some $(a',b')$ close to $(a,b)$, and we are back 
to the previous case. As $\|Q_0 - Q_{a,b}^{\omega,c}\|_X \le  
\|Q_0 - Q_{a,b}\|_X + \|Q_{a,b} - Q_{a,b}^{\omega,c}\|_X  
\le C\epsilon$, this shows that $\rho(Q(\cdot,t),Q_{a,b}^{\omega,c})  
\le C\epsilon$ for all $t \in \R$, and the result follows.  
This concludes the proof of Proposition~\ref{th:Qprop}.\QED 
 
\begin{Remark}\label{r:basis} 
In \cite{GSS2} the authors use the decomposition $X = \TT_{a,b} 
\oplus \tilde\NN_{a,b}$, where $\tilde \NN_{a,b}$ is the 
two-dimensional space spanned by 
\[ 
  \partial_\omega Q_{a,b} \,=\, \frac{\partial}{\partial \omega} 
  Q_{a,b}^{\omega,c}\Big|_{(\omega,c)=(0,0)}~, \quad 
  \partial_c Q_{a,b} \,=\, \frac{\partial}{\partial c} 
  Q_{a,b}^{\omega,c}\Big|_{(\omega,c)=(0,0)}~. 
\] 
This alternative decomposition has the advantage that $\langle H_{a,b}u,v 
\rangle = 0$ for all $u \in \tilde \NN_{a,b}$ and all $v \in \TT_{a,b}$, 
because 
\[ 
  H_{a,b}(\partial_\omega Q_{a,b}) \,=\, N'(Q_{a,b}) \,=\, Q_{a,b}~, \quad 
  H_{a,b}(\partial_c Q_{a,b}) \,=\, M'(Q_{a,b}) \,=\, \rmi 
  \partial_z Q_{a,b}~. 
\] 
Using in addition \reff{e:dNM} we also find 
\[ 
  \HH_{a,b} \,=\, -\pmatrix{ 
  \langle H_{a,b}(\partial_\omega Q_{a,b}),\partial_\omega Q_{a,b}\rangle & 
  \langle H_{a,b}(\partial_\omega Q_{a,b}),\partial_c Q_{a,b}\rangle \cr 
  \langle H_{a,b}(\partial_c Q_{a,b}),\partial_\omega Q_{a,b}\rangle & 
  \langle H_{a,b}(\partial_c Q_{a,b}),\partial_c Q_{a,b}\rangle}~. 
\] 
Remark that the spaces $\NN_{a,b}$ and $\tilde \NN_{a,b}$ are very close 
when $(a,b)$ is small, since 
\[ 
  \partial_\omega Q_{a,b} \,=\, \frac1{6ab}(a + b \rme^{-\rmi z})(1 + 
  \rmO(a^2+b^2))~, \quad 
  \partial_c Q_{a,b} \,=\, \frac1{6ab}(2a - b \rme^{-\rmi z})(1 + 
  \rmO(a^2+b^2))~. 
\] 
\end{Remark} 
 
\begin{Remark}[${2n\pi}$--periodic perturbations]\label{r:op} 
As an intermediate step between the periodic set-up considered in 
Theorem~\ref{th:orbit} and the case of arbitrary bounded 
perturbations for which no result is available so far, one can try 
to study the orbital stability of the travelling waves of 
\reff{e:nls} with respect to perturbations whose periods are 
integer multiples of the period of the original wave. This amounts 
to replacing the space $X$ in Proposition~\ref{th:Qprop} by 
$H^1_\per([0,2n\pi],\C)$ for some $n \ge 2$. In that case most of 
the results above remain valid, but the linear operator $H_{a,b}$ 
has now $2n-1$ negative eigenvalues. Thus the number of negative 
eigenvalues of $H_{a,b}$ minus the number of positive eigenvalues 
of the Hessian matrix $\HH_{a,b}$ is equal to $2n-2$, a nonzero 
even integer. This means that neither the Stability Theorem nor 
the Instability Theorem in \cite{GSS2} applies if $n \ge 2$. The 
only way out of this difficulty would be to replace the manifolds 
$\Sigma_{a,b}$ defined in \reff{e:Sabdef} by invariant manifolds 
higher codimension, which amounts to use additional conserved 
quantities of \reff{e:nls} instead of $N$ and $M$ only. 
\end{Remark} 
 
\section{Spectral stability} 
\label{s:ss} 
 
In this section we prove the spectral stability result in 
Theorem~\ref{th:spec}. We start with the evolution equation 
\reff{e:Qeq}, which we linearize about the stationary solution  
$Q_{a,b}(z)$ corresponding to the periodic wave $U_{a,b}(x,t) =  
\rme^{\rmi (p_{a,b}x-t)} Q_{a,b}(2k_{a,b}x)$ of the NLS equation. 
 We find the linear operator 
\begin{equation}\label{e:Aabdef} 
  \mathcal A_{a,b} Q = 4 \rmi k_{a,b}^2 Q_{zz} - 4p_{a,b}k_{a,b} 
  Q_z + \rmi (1-p_{a,b}^2)Q - 2\rmi |Q_{a,b}|^2Q - \rmi 
  Q_{a,b}^2 \overline Q~, 
\end{equation} 
which we consider in either the real Hilbert space $Y=L^2(\R,\C)$ 
(localized perturbations) or the real Banach space $Y=C_b(\R,\C)$ 
(bounded pertubations). To study the spectrum $\mathcal A_{a,b}$, it 
is convenient to decompose the elements of $Y$ into real and imaginary 
parts, in which case we obtain the matrix operator 
\begin{equation}\label{e:l_op} 
  \mathcal A_{a,b} = \left(\begin{array}{cc} 
  -4p_{a,b}k_{a,b}\partial_z+ 2R_{a,b}I_{a,b} 
  &-4k_{a,b}^2\partial_{zz} + (p_{a,b}^2-1)+ R_{a,b}^2 + 3I_{a,b}^2\\ 
  4k_{a,b}^2\partial_{zz} - (p_{a,b}^2-1) - 3R_{a,b}^2 - I_{a,b}^2 & 
  - 4p_{a,b}k_{a,b}\partial_z - 2R_{a,b}I_{a,b}\end{array}\right)~, 
\end{equation} 
where $Q_{a,b}=R_{a,b}+\rmi I_{a,b}$. We are now interested in the 
spectrum of this matrix operator in the (complexified) spaces 
$L^2(\R,\C^2)$ and $C_b(\R,\C^2)$. We prove that the 
spectrum of $\mathcal A_{a,b}$ in both spaces lies entirely on the 
imaginary axis, if $\|(a,b)\|$ is sufficiently small. This  
means that the periodic wave $U_{a,b}$ is spectrally stable  
in $Y$.  
 
\subsection{Bloch-wave decomposition and symmetries} 
\label{s:ss1} 
 
The spectral analysis of $\mathcal A_{a,b}$ relies upon the 
so-called Bloch-wave decomposition for differential operators with 
periodic coefficients. This method allows to show that the  
spectrum of $\mathcal A_{a,b}$ is exactly the same in both  
spaces $L^2(\R,\C^2)$ and $C_b(\R,\C^2)$, and can be described  
as the union of the point spectra of a family of operators with  
compact resolvent (see e.g. \cite{Mi,RS}). In our case, the  
operator $\mathcal A_{a,b}$ has $2\pi$-periodic coefficients and  
its spectrum in both $L^2(\R,\C^2)$ and $C_b(\R,\C^2)$ is given by 
\begin{equation}\label{e:bloch} 
  \sigma(\mathcal A_{a,b}) = \bigcup_{\gamma\in(-\frac12,\frac12]}\, 
  \sigma(\mathcal A_{a,b,\gamma})~, 
\end{equation} 
where the Bloch operators 
\[ 
  \mathcal A_{a,b,\gamma} = \left(\!\!\begin{array}{cc} 
  -4p_{a,b}k_{a,b}(\partial_z+\rmi\gamma)+ 2R_{a,b}I_{a,b} 
  &\!\!\!-4k_{a,b}^2(\partial_z+\rmi\gamma)^2 + (p_{a,b}^2-1)+ 
  R_{a,b}^2 + 3I_{a,b}^2\\ 
  4k_{a,b}^2(\partial_z+\rmi\gamma)^2 - (p_{a,b}^2-1) - 3R_{a,b}^2 - 
  I_{a,b}^2 &\!\!\! - 4p_{a,b}k_{a,b}(\partial_z+\rmi\gamma) - 
  2R_{a,b}I_{a,b}\end{array}\!\!\right) 
\] 
are linear operators in the Hilbert space of $2\pi$-periodic 
functions $L^2_\per([0,2\pi],\C^2)$. We can now reformulate the 
spectral result of Theorem~\ref{th:spec} as follows: 
 
\begin{Proposition}\label{p:abg} 
There exists $\delta_1>0$ such that, for any $\gamma \in  
(-\frac12,\frac12]$ and any $(a,b)\in\R^2$ with $\|(a,b)\| \le  
\delta_1$, the spectrum of the operator $\mathcal A\abg$ in  
$L^2_\per([0,2\pi],\C^2)$ satisfies $\sigma(\mathcal A\abg) 
\subset \rmi\R$. 
\end{Proposition} 
 
We equip the Hilbert space $L^2_\per([0,2\pi],\C^2)$ with the 
usual scalar product defined through 
\[ 
\langle (Q_1,Q_2)^t,(R_1,R_2)^t \rangle = \int_{0}^{2\pi} 
\Bigl(Q_1(z) 
  \overline R_1(z)+Q_2(z)\overline R_2(z)\Bigr) \dd z~. 
\] 
The operators $\mathcal A\abg$ are closed in this space with 
compactly embedded domain $H^2_\per([0,2\pi],\C^2)$. An immediate 
consequence of the latter property is that these operators have 
compact resolvent, so that their spectra are purely point spectra 
consisting of isolated eigenvalues with finite algebraic 
multiplicities. Our problem consists in locating these 
eigenvalues. 
 
The spectra of the operators $\mathcal A_{a,b}$ and $\mathcal 
A\abg$ possess several symmetries originating from the discrete 
symmetries of \reff{e:nls} and the symmetries of the wave profile 
$Q_{a,b}$. First, since the operator $\mathcal A_{a,b}$ has real 
coefficients, its spectrum is symmetric with respect to the real  
axis: $\sigma(\mathcal A_{a,b}) = \overline{\sigma(\mathcal A_{a,b})}$.   
For the Bloch operator $\mathcal A\abg$, the corresponding  
property is $\sigma(\mathcal A\abg) = \overline{\sigma 
(\mathcal A_{a,b,-\gamma})}$. Next, it is straightforward to  
check that $\mathcal A_{a,b}$ has a {\it reversibility symmetry},  
i.e. it anticommutes with the isometry $\mathcal S$ defined by 
\begin{equation}\label{e:s} 
  \mathcal S\pmatrix{Q_1(z) \cr Q_2(z)} \,=\,  
  \pmatrix{Q_1(-z) \cr -Q_2(-z)}~. 
\end{equation} 
Thus $\mathcal S\mathcal A_{a,b} = -\mathcal A_{a,b}\mathcal S$,  
which implies that the spectrum of $\mathcal A_{a,b}$ is  
symmetric with respect to the origin in the complex plane: 
$\sigma(\mathcal A_{a,b}) = -{\sigma(\mathcal A_{a,b})}$.  
The corresponding property for the Bloch operators is  
$\mathcal S\mathcal A\abg = -\mathcal A_{a,b,-\gamma}\mathcal S$,  
which implies that $\sigma(\mathcal A\abg) = -{\sigma(\mathcal  
A_{a,b,-\gamma})}$. Summarizing, the spectrum of $\mathcal 
A_{a,b}$ is symmetric with respect to both the real and the 
imaginary axis, and the spectra of the Bloch operators  
$\mathcal A\abg$ satisfy 
\begin{equation}\label{e:symg1} 
  \sigma(\mathcal A\abg) \,=\, \overline{\sigma(\mathcal 
  A_{a,b,-\gamma})} \,=\, -{\sigma(\mathcal A_{a,b,-\gamma})} \,=\, 
  -\overline{\sigma(\mathcal A_{a,b,\gamma})}~. 
\end{equation} 
In particular, the spectrum of $\mathcal A\abg$ is symmetric with 
respect to the imaginary axis and we can restrict ourselves to 
positive values $\gamma\in[0,\frac12]$ without loss of generality. 
 
Using now the relations \reff{e:symQ} for the wave profile 
$Q_{a,b}$, we see that the spectra of $\mathcal A_{a,b}$ satisfy 
$\sigma(\mathcal A_{a,b}) = \sigma(\mathcal A_{-a,b}) = 
\sigma(\mathcal A_{-a,-b})$ and $\sigma(\mathcal A_{b,a}) = - 
\overline{\sigma(\mathcal A_{a,b})}$. (Actually, the last  
equality is easier to establish if we use the complex form  
\reff{e:Aabdef} of the operator $\mathcal A_{a,b}$, for which  
we have $\mathcal A_{b,a}(\rme^{-\rmi z}Q) = - \rme^{-\rmi z}  
\overline{\mathcal A_{a,b}}Q$.) Similarly, we find for the Bloch  
operators 
\begin{equation}\label{e:symg2} 
  \sigma(\mathcal A_{a,b,\gamma}) \,=\, \sigma(\mathcal A_{-a,b,\gamma}) 
  \,=\, \sigma(\mathcal A_{-a,-b,\gamma})~,\quad \mbox{and}\quad 
  \sigma(\mathcal A_{b,a,\gamma}) \,=\, -\overline{\sigma(\mathcal 
  A_{a,b,\gamma})} ~. 
\end{equation} 
 
Finally, we note the formal relation $\mathcal A_{a,b} = -\rmi 
H_{a,b}$ between the linearized operator \reff{e:Aabdef} and the 
second variation of the energy defined in \reff{e:Habdef}.  
For the Bloch operators, we can write in a similar way 
\[ 
  \mathcal A_{a,b,\gamma} \,=\, JH_{a,b,\gamma}~,\quad J \,=\, 
  \pmatrix{0 & 1 \cr -1 & 0}~, 
\] 
with 
\[ 
  H_{a,b,\gamma} =\! \left(\!\!\begin{array}{cc} 
  \!-4k_{a,b}^2(\partial_z+\rmi\gamma)^2 + (p_{a,b}^2-1) + 
  3R_{a,b}^2 + I_{a,b}^2 & 
  \!\!\!\!4p_{a,b}k_{a,b}(\partial_z+\rmi\gamma) 
  + 2R_{a,b}I_{a,b}\\ 
  \!-4p_{a,b}k_{a,b}(\partial_z+\rmi\gamma)+ 2R_{a,b}I_{a,b} 
  &\!\!\!\!-4k_{a,b}^2(\partial_z+\rmi\gamma)^2 + (p_{a,b}^2-1)+ 
  R_{a,b}^2 + 3I_{a,b}^2\end{array}\!\!\right)\!. 
\] 
Actually, this property is a consequence of the Hamiltonian 
structure of the NLS equation. Though some properties induced by 
this structure are exploited, we shall not make an explicit use of 
the Hamiltonian structure itself in the proof of spectral stability.

\subsection{First perturbation argument and properties of the unperturbed operators} 
 
Our spectral analysis for the operators $\mathcal A_{a,b,\gamma}$ 
relies upon perturbation arguments in which we regard $\mathcal 
A_{a,b,\gamma}$ as small bounded perturbations of the operators 
with constant coefficients 
\[ 
  \mathcal A^0_{a,b,\gamma} \,=\, \left(\begin{array}{cc} 
  -4p_{a,b}k_{a,b}(\partial_z+\rmi\gamma) 
  &-4k_{a,b}^2(\partial_z+\rmi\gamma)^2 + (p_{a,b}^2-1)\\ 
  4k_{a,b}^2(\partial_z+\rmi\gamma)^2 - (p_{a,b}^2-1) & - 
  4p_{a,b}k_{a,b}(\partial_z+\rmi\gamma) \end{array}\right)~. 
\] 
The difference $\mathcal A^1_{a,b}:=\mathcal 
A_{a,b,\gamma}-\mathcal A^0_{a,b,\gamma}$ is a bounded operator 
with norm $\|\mathcal A^1_{a,b}\|=\rmO(a^2+b^2)$, as 
$(a,b)\to(0,0)$. 
 
A straightforward Fourier analysis allows to compute the spectra 
of the operators $A^0_{a,b,\gamma}$: 
\begin{equation}\label{e:spec0} 
  \sigma(\mathcal A\abg^0) = \{ \rmi\omega_{a,b,\gamma}^{\pm,n}, \ 
  \omega_{a,b,\gamma}^{\pm,n}= -4p_{a,b}k_{a,b}(n+\gamma) \pm \left( 
  4k_{a,b}^2(n+\gamma)^2 + p_{a,b}^2-1\right) ,\ n\in \Z \} \subset 
  \rmi\R~, 
\end{equation} 
in which the eigenvalues are all semi-simple with eigenfunctions 
\[ 
  e^{\pm,n} \,=\, \rme^{\rmi nz}\pmatrix{1 \cr \pm\rmi}~, \quad  
  \mathcal A^0_{a,b,\gamma}e^{\pm,n} \,=\, 
  \rmi\omega_{a,b,\gamma}^{\pm,n} \, e^{\pm,n}~. 
\] 
Furthermore, the resolvent operators $\RR^0_{a,b,\gamma}(\lambda) = 
(\lambda\1 - \mathcal A^0_{a,b,\gamma})^{-1}$ have norms 
\[ 
  \|\RR\ze(\lambda)\| \,=\, \frac1{{\rm dist}(\lambda,\sigma(\mathcal 
  A\abg^0))}~, \quad \lambda \notin \sigma(\mathcal A\abg^0)~. 
\] 
A simple perturbation argument shows now that the spectrum of 
$\mathcal A_{a,b,\gamma}$ stays close to $\sigma(\mathcal 
A\abg^0)$ provided $\|(a,b)\|$ is sufficiently small. More 
precisely, we have the following result. 
 
\begin{Lemma}\label{l:pert} 
For any $c>0$ there exists $\delta>0$ such that for any 
$\gamma\in[0,\frac12]$ and any $(a,b) \in \R^2$ with 
$\|(a,b)\|\leq\delta$ the spectrum of $\mathcal A_{a,b,\gamma}$ 
satisfies 
\[ 
  \sigma(\mathcal A\abg) \,\subset\, \bigcup_{n\in\Z} 
  B(\rmi\omega_{a,b,\gamma}^{-,n};c)\,\cup\, \bigcup_{n\in\Z} 
  B(\rmi\omega_{a,b,\gamma}^{+,n};c), 
\] 
in which $B(\rmi\omega_{a,b,\gamma}^{\pm,n};c)$ represents the 
open ball centered at $\rmi\omega_{a,b,\gamma}^{\pm,n}$ with 
radius $c$. 
\end{Lemma} 
 
\begin{Proof} 
For any $\lambda\notin \bigcup_{n\in\Z} B(\rmi\omega_{a,b, 
\gamma}^{\pm,n};c)$, we write 
\[ 
  \lambda\1 - \mathcal A\abg \,=\, \lambda\1 - \mathcal A\ze - \mathcal 
  A^1_{a,b} \,=\, (\lambda\1 - \mathcal A\ze) \left(\1 - \RR\ze(\lambda) 
  \mathcal A^1_{a,b}\right)~. 
\] 
Since 
\[ 
  \|\RR\ze(\lambda) \mathcal A^1_{a,b}\|\leq \frac1c \,\| \mathcal 
  A^1_{a,b}\|~,\quad \| \mathcal A^1_{a,b}\| \,=\,\rmO(a^2+b^2)~, 
\] 
upon choosing $\delta$ sufficiently small, we have that $\1 - 
\RR\ze(\lambda) \mathcal A^1_{a,b}$ is invertible, so that 
$\lambda\1 - \mathcal A\abg$ is invertible, as well. This proves 
that $\lambda$ does not belong to $\sigma(\mathcal A\abg)$. 
\end{Proof} 
 
In order to locate the spectra of $\mathcal A\abg$, we need a more 
precise description of the spectra of $\mathcal A\ze$. Looking at 
$a=b=0$ we find that: 
\begin{itemize} 
\item if $\gamma=0$, all nonzero eigenvalues of 
  $\mathcal A^0_{0,0,0}$ are double, 
  \[ 
    \omega_{0,0,0}^{+,n} \,=\, \omega_{0,0,0}^{+,1-n}~,\quad 
    \omega_{0,0,0}^{-,n} \,=\, \omega_{0,0,0}^{-,-1-n}~, 
  \] 
  and zero is an eigenvalue of multiplicity 4, 
  \[ 
    \omega_{0,0,0}^{\pm,0} \,=\, \omega_{0,0,0}^{+,1} \,=\, 
    \omega_{0,0,0}^{-,-1} \,=\, 0~; 
  \] 
  \item if $0<\gamma<\frac12$ all eigenvalues are simple; 
  \item if $\gamma=\frac12$, there is a pair of simple 
  eigenvalues $\pm\rmi$, 
  \[ 
    \omega_{0,0,\frac12}^{-,-1} \,=\, - \omega_{0,0,\frac12}^{+,0} \,=\, 
    1~, 
  \] 
  and the other eigenvalues are all double, 
    \[ 
    \omega_{0,0,\frac12}^{+,n} \,=\, \omega_{0,0,\frac12}^{+,-n}~,\quad 
    \omega_{0,0,\frac12}^{-,n} \,=\, \omega_{0,0,\frac12}^{-,-2-n}~. 
  \] 
\end{itemize} 
We therefore distinguish three cases: $\gamma\approx 0$, 
$\gamma\approx \frac12$, and $\gamma \in [\gamma_*,\frac12- \gamma_*]$ 
for some $\gamma_* \in (0,\frac14)$, which we treat separately in the 
next paragraphs.  In each case, the starting point is an estimate of 
the distance between any pair of eigenvalues of $\mathcal A\ze$, which 
is directly obtained from the explicit formulas \reff{e:spec0}.  We 
use this estimate to construct an infinite family of mutually disjoint 
sets (balls or finite unions of balls) with the property that the 
spectrum of $\mathcal A\abg$ is contained in their union.  Inside each 
set $\mathcal A\abg$ will have a finite number of eigenvalues (one, 
two or four) so that the problem reduces to showing that these 
eigenvalues are purely imaginary. In Propositions~\ref{p:simple}, 
\ref{p:far}, \ref{p:zero}, and \ref{p:12} below, we show that in all 
three cases the spectrum of $\mathcal A\abg$ is purely imaginary, 
provided $\|(a,b)\|$ is sufficiently small. This proves 
Proposition~\ref{p:abg}.

\subsection{Spectrum for $\gamma$ away from $0$ and $\frac12$} 
 
We start with the case $\gamma\in[\gamma_*,\frac12-\gamma_*]$, 
when the operators with constant coefficients $\mathcal A\abg^0$ 
have only simple eigenvalues: 
 
\begin{Lemma}\label{l:p1} 
For any $\gamma_*\in(0,\frac14)$, there exist positive constants 
$c_*$ and $\delta_*$ such that for any 
$\gamma\in[\gamma_*,\frac12-\gamma_*]$ and any $(a,b)$ with 
$\|(a,b)\|\leq\delta_*$, we have 
\[ 
|\rmi \omega\abg^{\sigma,n} - \rmi \omega\abg^{\tau,p}| \geq 
c_*~,\quad \hbox{for all} \ p,n\in\Z \ \mbox{and all}\ 
\sigma,\tau\in\{-,+\}\ \mbox{with}\ (\sigma,n)\not=(\tau,p)~. 
\] 
\end{Lemma} 
 
This lemma shows that the distance between any pair of eigenvalues 
of $\mathcal A\ze$ is strictly positive, uniformly for small 
$\|(a,b)\|$ and $\gamma\in[\gamma_*,\frac12-\gamma_*]$. This allows us 
to find an infinite sequence of mutually disjoint balls with the 
property that the spectrum of $\mathcal A\abg$ is contained in their 
union, and that inside each ball both operators have precisely one 
simple eigenvalue. The symmetry of the spectrum of $\mathcal A\abg$ 
with respect to the imaginary axis then implies that this simple 
eigenvalue is purely imaginary. We point out that classical 
perturbation results for families of simple eigenvalues \cite{K1} do 
not directly apply, since here we have infinitely many eigenvalues. 
 
\begin{Proposition}\label{p:simple} 
Fix $\gamma_*\in(0,\frac14)$. Then there exist positive  constants 
$c$ and $\delta$ such that for any 
$\gamma\in[\gamma_*,\frac12-\gamma_*]$ and any $(a,b)$ with 
$\|(a,b)\|\leq\delta$, the following properties hold: 
\begin{enumerate} 
\item The spectrum  of $\mathcal A\abg$ satisfies 
\[ 
\sigma(\mathcal A\abg) \subset \bigcup_{n\in\Z} 
B(\rmi\omega\abg^{-,n};c)\,\cup\, \bigcup_{n\in\Z} 
B(\rmi\omega\abg^{+,n};c) ~, 
\] 
and the closed balls $\overline{B(\rmi\omega\abg^{\pm,n};c)}$ are 
mutually disjoints. 
 
\item Inside each ball $B(\rmi\omega\abg^{\pm,n};c)$ the operator 
$\mathcal A\abg$ has precisely one eigenvalue, which is purely 
imaginary.\footnote{Here and in the rest of the paper we say that 
``$\mathcal A\abg$ has $n$ eigenvalues inside the set $B$'' when 
the sum of the algebraic multiplicities of the eigenvalues of 
$\mathcal A\abg$ inside $B$ is equal to $n$.} 
\end{enumerate} 
\end{Proposition}

\begin{Proof} 
(i) We may choose any $c \leq c_*/4$ with $c_*$ the constant in 
Lemma~\ref{l:p1}, so that the balls 
$\overline{B(\rmi\omega\abg^{\pm,n};c)}$ are mutually disjoints, 
and then apply Lemma~\ref{l:pert}. 
 
(ii) Consider a ball $B(\rmi\omega\abg^{\pm,n};c)$. Inside this 
ball $\mathcal A\ze$ has precisely one eigenvalue 
$\rmi\omega\abg^{\pm,n}$ with associated spectral projection 
$\Pi^{0,n}\abg$ satisfying $\|\Pi^{0,n}\abg\|=1$. The result (i) 
provides us with a spectral decomposition for the operator 
$\mathcal A\abg$, and we can compute the spectral projection 
$\Pi^n\abg$ associated to $B(\rmi\omega\abg^{\pm,n};c)$ via the 
Dunford integral formula 
\begin{equation}\label{e:dun} 
  \Pi^n\abg \,=\, \frac1{2\pi\rmi} \oint_{C_n} \RR\abg(\lambda)\mathrm 
  d\lambda~, 
\end{equation} 
in which $C_n$ is the boundary of $B(\rmi\omega\abg^{\pm,n};c)$ 
and $\RR\abg(\lambda)=(\lambda\1- \mathcal A\abg)^{-1}$. Using the 
formula for the resolvent 
\[ 
  \RR\abg(\lambda) \,=\, \RR\ze(\lambda)\left(\1 - \mathcal A^1_{a,b} 
  \RR\ze(\lambda)\right)^{-1}~, 
\] 
which holds for sufficiently small $\delta$ since $\|\mathcal 
A^1_{a,b}\|=\rmO(a^2+b^2)$, we compute the difference 
\[ 
  \Pi^n\abg-\Pi^{0,n}\abg \,=\, \frac1{2\pi\rmi} \oint_{C_n} 
  \RR\ze(\lambda)\sum_{k\geq1} \left( \mathcal A^1_{a,b} \RR\ze 
  (\lambda)\right)^k\mathrm d\lambda~. 
\] 
Since $\|\RR\ze(\lambda)\|=1/c$, for $\lambda\in C_n$, we conclude 
that 
\[ 
  \|\Pi^n\abg-\Pi^{0,n}\abg\|\, \leq\, \sum_{k\geq1} \left( 
  \frac1c\, \|\mathcal A^1_{a,b}\| \right)^k \,=\, \frac{\|\mathcal 
  A^1_{a,b}\|}{c-\|\mathcal A^1_{a,b}\|}~. 
\] 
Upon choosing $\delta$ small enough we achieve 
\begin{eqnarray*} 
  \|\Pi^n\abg-\Pi^{0,n}\abg\|& < & \frac1{1+ 
  \|\Pi^n\abg-\Pi^{0,n}\abg\|} \,=\, \frac1{\|\Pi^{0,n}\abg\|+ 
  \|\Pi^n\abg-\Pi^{0,n}\abg\|}\\ & \leq & \min\left( 
  \frac1{\|\Pi^{0,n}\abg\|}, \frac1{\|\Pi^{n}\abg\|} \right)~, 
\end{eqnarray*} 
so that the projections $\Pi^n\abg$ and $\Pi^{0,n}\abg$ realize 
isomorphisms between the associated spectral subspaces of $\mathcal 
A\ze$ and $\mathcal A\abg$ (\cite[Lemma B.1]{HLS}; see also 
\cite[Chapter I {\S}6.8]{K1}). In particular, they have the same 
finite rank, so that $\mathcal A\abg$ has precisely one simple 
eigenvalue inside the ball $B(\rmi\omega\abg^{\pm,n};c)$. Finally, 
since the spectrum is symmetric with respect to the imaginary axis 
(\ref{e:symg1}) this simple eigenvalue is necessarily purely 
imaginary, which concludes the proof. 
\end{Proof} 
 
\subsection{Spectrum for small $\gamma$}\label{ss:small} 
 
We start again by analyzing the distance between the eigenvalues 
of $\mathcal A\ze$, now for small values of $\gamma$. Since at 
$a=b=\gamma=0$ the spectrum of $\mathcal A\ze$ consists of double 
nonzero eigenvalues and a quadruple eigenvalue at zero, for small 
$a$, $b$, and $\gamma$ we expect pairs of arbitrarily close 
eigenvalues together with four eigenvalues close to the origin. A 
precise description of the location of these eigenvalues is given 
in the following lemma. 
 
\begin{Lemma}\label{l:p2} 
There exist positive constants $\gamma_0$, $c$, and $\delta$, such 
that the following holds, for any $\gamma\in[0,\gamma_0]$, and 
$(a,b)$ with $\|(a,b)\|\leq\delta$: 
\begin{enumerate} 
  \item $\rmi\omega\abg^{\pm,0},\ \rmi\omega\abg^{+,1}, \ 
  \rmi\omega\abg^{-,-1}\in B(0;1)$; 
  \item $\rmi\omega\abg^{-,1},\ \rmi\omega\abg^{+,-1}, \ 
  \rmi\omega\abg^{\pm,n}\notin B(0;4)$, $n\in\Z\setminus 
  \{-1,0,1\}$; 
  \item $|\rmi\omega\abg^{+,n} - \rmi\omega\abg^{+,p}|\geq c$, 
  $n,p\in \Z\setminus\{0,1\}$, $p\neq n$, $p\neq 1-n$; 
  \item $|\rmi\omega\abg^{-,n} - \rmi\omega\abg^{-,p}|\geq c$, 
  $n,p\in \Z\setminus\{-1,0\}$, $p\neq n$, $p\neq -1-n$; 
  \item $|\rmi\omega\abg^{+,n} - \rmi\omega\abg^{-,p}|\geq c$, 
  $n\in \Z\setminus\{0,1\},\ p\in\Z\setminus\{-1,0\}$~. 
\end{enumerate} 
\end{Lemma} 
 
The first two properties (i)-(ii) in this lemma together with the 
perturbation result in Lemma~\ref{l:pert} provides us with a 
spectral splitting for $\mathcal A\abg$: 
\[ 
  \sigma(\mathcal A\abg) \,=\, \sigma_1(\mathcal A\abg) \,\cup\, 
  \sigma_2(\mathcal A\abg)~, 
\] 
with 
\[ 
  \sigma_1(\mathcal A\abg) \,\subset\, B(0;2)~,\quad \sigma_2(\mathcal 
  A\abg) \cap B(0;3) \,=\, \emptyset~. 
\] 
Inside the ball $B(0;2)$ we find the part of the spectrum of $\mathcal 
A\abg$ which is close to the quadruple zero eigenvalue of $\mathcal 
A_{0,0,0}^0$, whereas the rest of the spectrum lies outside the ball 
$B(0;3)$. The last properties (iii)-(v) show that the eigenvalues 
outside $B(0;3)$ are well separated except for the pairs 
$(\rmi\omega\abg^{+,n},\rmi\omega\abg^{+,1-n})$, $n\in 
\Z\setminus\{0,1\}$, and $(\rmi\omega\abg^{-,n},\rmi 
\omega\abg^{-,-1-n})$, $n\in \Z\setminus\{-1,0\}$, which may be 
arbitrarily close. At $a=b=\gamma=0$, these are precisely the double 
eigenvalues of $\mathcal A^0_{0,0,0}$. Notice however that, for fixed 
$a$, $b$, $\gamma$, the distances 
$|\rmi\omega\abg^{+,n}-\rmi\omega\abg^{+,1-n}|$ and 
$|\rmi\omega\abg^{-,n}-\rmi\omega\abg^{-,-1-n}|$ typically grow like 
$\rmO(n^2)$ as $|n| \to \infty$, see \reff{e:spec0}. 
 
We analyze these two parts $\sigma_1(\mathcal A\abg)$ and 
$\sigma_2(\mathcal A\abg)$ of the spectrum of $\mathcal A\abg$ 
separately in the Propositions~\ref{p:far} and \ref{p:zero} below. 
 
\begin{Proposition}\label{p:far} 
There exist positive constants $\gamma_0$, $c$, and $\delta$, such 
that for any $\gamma\in[0,\gamma_0]$, and $(a,b)$ with 
$\|(a,b)\|\leq\delta$, the following holds. 
\begin{enumerate} 
\item The spectrum $\sigma_2(\mathcal A\abg)$ satisfies 
\[ 
  \sigma_2(\mathcal A\abg)\, \subset\, 
  B(\rmi\omega\abg^{-,1};c)\,\cup\, B(\rmi\omega\abg^{+,-1};c)\, 
  \cup \bigcup_{n\not=\pm1,0} B(\rmi\omega\abg^{-,n};c)\, \cup\, 
  \bigcup_{n\not=\pm1,0} B(\rmi\omega\abg^{+,n};c)~, 
\] 
in which the balls $B(\rmi\omega\abg^{\pm,n};c)$ are mutually 
disjoints, except for some pairs 
$(\rmi\omega\abg^{+,n},\rmi\omega\abg^{+,1-n})$, $n\in 
\Z\setminus\{0,1\}$, or 
$(\rmi\omega\abg^{-,n},\rmi\omega\abg^{-,-1-n})$, $n\in 
\Z\setminus\{-1,0\}$. 
 
\item Inside each ball $B(\rmi\omega\abg^{\pm,n};c)$ the operator 
$\mathcal A\abg$ has either one or two eigenvalues, which are 
purely imaginary. 
\end{enumerate} 
\end{Proposition} 
 
\begin{Proof} 
The result (i) is obtained from Lemma~\ref{l:p2}~(iii)-(v) and 
Lemma~\ref{l:pert}, just as the first part of 
Proposition~\ref{p:simple}. The only difference is that here we 
have pairs of balls which are not disjoint. As was noticed above,  
the distances $|\rmi\omega\abg^{+,n}-\rmi\omega\abg^{+,1-n}|$ and 
$|\rmi\omega\abg^{-,n}-\rmi\omega\abg^{-,-1-n}|$ grow like $\rmO(n^2)$  
as $|n| \to \infty$, so that we have in general a finite number  
of such pairs for a given value of $a,b,\gamma$.  
 
(ii) For the balls $B(\rmi\omega\abg^{\pm,n};c)$ which are 
disjoint from all the others we can argue and conclude as in the proof 
of Proposition~\ref{p:simple}. It remains to consider the case of 
two balls which are not disjoint. Choose a pair of eigenvalues 
$(\rmi\omega\abg^{+,n},\rmi\omega\abg^{+,1-n})$ such that 
$B(\rmi\omega\abg^{+,n};c) \cap B(\rmi\omega\abg^{+,1-n};c)  
\neq \emptyset$ (the argument is similar for a pair 
$(\rmi\omega\abg^{-,n},\rmi\omega\abg^{-,-1-n})$). We construct 
the spectral projection $\Pi^{n,1-n}\abg$ for $\mathcal A\abg$ 
corresponding to the union of these balls with the help of the 
Dunford integral formula \reff {e:dun}, in which the circle $C_n$  
is replaced by the smallest circle $C_r$ with radius 
$c<r<2c$, centered on the imaginary axis, which contains  
both balls. The spectral projection $\Pi^{0,n,1-n}\abg$ for 
$\mathcal A\ze$ has unit norm again, and since 
$\|\RR\ze(\lambda)\|\leq 1/c$ for $\lambda\in C_r$, we  
easily find 
\begin{equation}\label{e:pi} 
  \|\Pi^{n,1-n}\abg-\Pi^{0,n,1-n}\abg\| \,\leq\, \frac 
  rc\,\frac{\|\mathcal A^1_{a,b}\|}{c-\|\mathcal A^1_{a,b}\|} \,\leq\, 
  \frac{2\|\mathcal A^1_{a,b}\|}{c-\|\mathcal A^1_{a,b}\|} \,=\, 
  \rmO(a^2+b^2)~. 
\end{equation} 
As in the proof of Proposition~\ref{p:simple} we now choose 
$\delta$ sufficiently small such that these projections realize 
isomorphisms between the associated spectral subspaces of 
$\mathcal A\ze$ and $\mathcal A\abg$. In particular, these 
subspaces have the same finite rank equal to 2, which proves that 
$\mathcal A\abg$ has precisely two eigenvalues in 
$B(\rmi\omega\abg^{+,n};c)\cup B(\rmi\omega\abg^{+,1-n};c)$. 
 
In order to show that these two eigenvalues do not move off the 
imaginary axis we choose an appropriate basis of the associated 
two-dimensional eigenspace and compute the $2\times2$ matrix 
representing the action of $\mathcal A\abg$ on this space. Then it 
suffices to show that this matrix has purely imaginary 
eigenvalues. We start with the basis 
\[ 
  \xi_{00} \,=\, \frac1{2\sqrt\pi}\, \rme^{\rmi nz} 
  \pmatrix{1 \cr \rmi}~,\quad \xi_{01} \,=\, \frac1{2\sqrt\pi}\,  
  \rme^{\rmi(1- n)z} \pmatrix{1 \cr \rmi}~, 
\] 
of the two-dimensional eigenspace of $\mathcal A\ze$, which 
satisfies $\langle-\rmi J\xi_{0k},\xi_{0\ell}\rangle = 
\delta_{k\ell}$. We claim that for $\mathcal A\abg$ we can find a 
basis with the same property. Indeed, consider the vectors 
\[ 
  \widetilde \xi_0 \,=\, \Pi^{n,1-n}\abg \xi_{00}~, \quad \widetilde 
  \xi_1 \,=\, \Pi^{n,1-n}\abg \xi_{01}~, 
\] 
which form a basis of the two-dimensional eigenspace of $\mathcal 
A\abg$. From \reff{e:pi} we obtain $\langle-\rmi J \widetilde 
\xi_{0},\widetilde\xi_{0}\rangle = 1 + \rmO(a^2+b^2) >0$, so that 
the vector $\xi_0$ defined by 
\[ 
  \xi_0 \,=\, \frac 1{\langle-\rmi 
  J\widetilde\xi_{0},\widetilde\xi_{0}\rangle^{1/2}} 
  \,\widetilde\xi_0~, 
\] 
satisfies $\langle-\rmi J\xi_{0},\xi_{0}\rangle = 1$. Then we define 
successively 
\[ 
  \widehat\xi_1 \,=\, \widetilde\xi_1 -\overline{\langle-\rmi 
  J\xi_{0},\widetilde\xi_{1}\rangle}\,\xi_0~,\quad \xi_1 \,=\, \frac 
  1{\langle-\rmi J\widehat\xi_{1},\widehat\xi_{1}\rangle^{1/2}} 
  \,\widehat\xi_1~, 
\] 
and find $\langle-\rmi J\xi_{0},\xi_{1}\rangle = 0$ and 
$\langle-\rmi J\xi_{1},\xi_{1}\rangle = 1$, which proves the 
claim. 
 
The property $\langle-\rmi J\xi_{k},\xi_{\ell}\rangle = 
\delta_{k\ell}$ implies that the action of $\mathcal A_{a,b,\gamma}$ 
on the two-dimensional space spanned by $\{\xi_0,\xi_1\}$ is given 
by the matrix 
\[ 
  \mathcal M_{a,b,\gamma} \,=\, \pmatrix{\langle \mathcal 
  A_{a,b,\gamma}\xi_0,-\rmi J \xi_0\rangle & \langle \mathcal 
  A_{a,b,\gamma}\xi_1,-\rmi J \xi_0\rangle \cr \langle \mathcal 
  A_{a,b,\gamma}\xi_0,-\rmi J \xi_1\rangle & \langle \mathcal 
  A_{a,b,\gamma}\xi_1,-\rmi J \xi_1\rangle}~. 
\] 
Using the decomposition $\mathcal A_{a,b,\gamma} = J 
H_{a,b,\gamma}$ we find 
\[ 
  \langle \mathcal A_{a,b,\gamma}\xi_k,-\rmi J \xi_l\rangle \,=\, 
  \langle JH_{a,b,\gamma}\xi_k,-\rmi J \xi_l\rangle \,=\,\langle 
  H_{a,b,\gamma}\xi_k,-\rmi J^{-1} J \xi_l\rangle \,=\, \rmi \langle 
  H_{a,b,\gamma}\xi_k,\xi_l\rangle ~, 
\] 
so that 
\[ 
  \mathcal M_{a,b,\gamma} \,=\, \rmi\pmatrix{\langle 
  H_{a,b,\gamma}\xi_0, \xi_0\rangle & \langle H_{a,b,\gamma}\xi_1, 
  \xi_0\rangle\cr \langle H_{a,b,\gamma}\xi_0, \xi_1\rangle & \langle 
  H_{a,b,\gamma}\xi_1, \xi_1\rangle}~. 
\] 
Since $\langle H_{a,b,\gamma} Q,R\rangle =  \overline{\langle 
H_{a,b,\gamma} R,Q\rangle}$, we conclude that this matrix always 
has purely imaginary eigenvalues. This completes the proof. 
\end{Proof} 
 
\begin{Remark}\label{r:krein} 
The last part of this proof is a simple version of the well-known 
result for general Hamiltonian systems which asserts that 
colliding purely imaginary eigenvalues do not leave the imaginary 
axis when they have the same Krein signature (see e.g. 
\cite{KKS}). In the case of the four eigenvalues close 
to the origin, which we treat in the next proposition, the same 
argument does not work anymore (these eigenvalues have opposite 
Krein signatures). Instead, we compute an explicit expansion of 
the restriction of $\mathcal A\abg$ to the associated eigenspace 
which allows to show that these four eigenvalues are purely 
imaginary. 
\end{Remark} 
 
\begin{Proposition}\label{p:zero} 
There exist positive constants $\gamma_0$, $c$, and $\delta$, such 
that for any $\gamma\in[0,\gamma_0]$, and $(a,b)$ with 
$\|(a,b)\|\leq\delta$, the set $\sigma_1(\mathcal A\abg)$ consists 
of four purely imaginary eigenvalues. 
\end{Proposition} 
 
\begin{Proof} 
As in the previous cases, upon choosing $\delta$ sufficiently 
small, we obtain that $\mathcal A\abg$ has precisely four 
eigenvalues inside the ball $B(0;2)$. In order to locate these 
four eigenvalues we construct a suitable basis for the associated 
eigenspace and compute the $4\times4$ matrix $\mathcal M\abg$ 
representing the action of $\mathcal A\abg$ on this space. Then we 
show that this matrix has purely imaginary eigenvalues. 
 
We start with the particular cases $a=b=0$, and $\gamma=0$. In the 
first case, the operator $\mathcal A_{0,0,\gamma} = \mathcal 
A^0_{0,0,\gamma}$ has constant coefficients, so that we can 
explicitly compute the basis and the matrix. We choose the real 
basis 
\[ 
  \xi^{(0)}_{0,0,\gamma} \,=\, \pmatrix{\sin z \cr \cos z}~, \quad 
  \xi^{(1)}_{0,0,\gamma} \,=\, \pmatrix{0 \cr 1}~, \quad  
  \xi^{(2)}_{0,0,\gamma} \,=\, \pmatrix{\cos z \cr -\sin z}~, \quad  
  \xi^{(3)}_{0,0,\gamma} \,=\, \pmatrix{1 \cr 0}~, 
\] 
in which we find 
\[ 
  \mathcal M_{0,0,\gamma} \,=\, \pmatrix{4\rmi\gamma\, 
  \mathbf D_2 & - 4\gamma^2 \mathbf 1_2 \cr 4\gamma^2 \mathbf 1_2 & 
  4\rmi\gamma\, \mathbf D_2}~,\quad  
  \mathbf D_2 \,=\, \pmatrix{1 & 0 \cr 0 & -1}~, \quad 
  \mathbf 1_2 \,=\, \pmatrix{1 & 0 \cr 0 & 1}~. 
\] 
 
Next, we consider the operator $\mathcal A_{a,b,0}$. As for the 
operator $H_{a,b}$ in Section~\ref{s:os}, we have that 
 \[ 
  \frac{\partial}{\partial \phi} \RR_{(\phi,\xi)} 
  Q_{a,b}\Big|_{(\phi,\xi)=(0,0)} \,=\, - \rmi Q_{a,b} ~, \quad 
  \frac{\partial}{\partial \xi} \RR_{(\phi,\xi)} 
  Q_{a,b}\Big|_{(\phi,\xi)=(0,0)}\,=\, \partial_z Q_{a,b} ~, 
\] 
belong to the kernel of $\mathcal A_{a,b,0}$. In addition, since 
\[ 
  H_{a,b}(\partial_\omega Q_{a,b}) \,=\, N'(Q_{a,b}) \,=\, Q_{a,b}~, \quad 
  H_{a,b}(\partial_c Q_{a,b}) \,=\, M'(Q_{a,b}) \,=\, \rmi 
  \partial_z Q_{a,b}~, 
\] 
(see Remark~\ref{r:basis}) and $\mathcal A_{a,b,0}=-\rmi H_{a,b}$, 
we have 
\[ 
  \mathcal A_{a,b,0}(\partial_\omega Q_{a,b}) \,=\, -\rmi Q_{a,b}~, \quad 
  \mathcal A_{a,b,0}(\partial_c Q_{a,b}) \,=\, 
  \partial_z Q_{a,b}~, 
\] 
which provides us with two principal vectors in the generalized 
kernel of $\mathcal A_{a,b,0}$. Together with the two vectors in 
the kernel of $\mathcal A_{a,b,0}$ these give us a basis for the 
four-dimensional  eigenspace of $\mathcal A_{a,b,0}$. At $a=b=0$ 
we must find the basis above so that we set 
\begin{eqnarray*} 
  && \xi^{(0)}_{a,b,0} \,=\, -\frac1a \partial_z Q_{a,b} \,=\, 
  \pmatrix{\sin z \cr \cos z} + \rmO(|b|(|a|+|b|))~, \\ 
  && \xi^{(1)}_{a,b,0} \,=\, \frac1b\left( \rmi Q_{a,b} + \partial_z 
  Q_{a,b}\right) \,=\, \pmatrix{0 \cr 1} + \rmO(|a|(|a|+|b|))~,\\ 
  && \xi^{(2)}_{a,b,0} \,=\, 2a (2\partial_\omega Q_{a,b} - \partial_c 
  Q_{a,b}) \,=\, \pmatrix{\cos z \cr -\sin z} + \rmO(a^2+b^2)~,\\ 
  && \xi^{(3)}_{a,b,0} \,=\, 2b (\partial_\omega Q_{a,b} + \partial_c 
  Q_{a,b}) \,=\, \pmatrix{1 \cr 0} + \rmO(a^2+b^2)~, 
\end{eqnarray*} 
and a straightforward calculation gives the matrix 
\[ 
  \mathcal M_{a,b,0} \,=\, \pmatrix{\mathbf 0_2 & \mathbf 
  M_2(a,b) \cr \mathbf 0_2 & \mathbf 0_2}~,\quad   
  \mathbf 0_2 \,=\, \pmatrix{0 & 0 \cr 0 & 0}~,\quad  
  \mathbf M_2 \,=\, \pmatrix{- 2a^2 & - 4ab \cr - 4ab & -2b^2} 
  + \rmO(a^4+b^4)~. 
\] 
 
Finally, we consider the full operator $\mathcal A_{a,b,\gamma}$ 
and construct a basis $\{\xi^{(0)}_{a,b,\gamma}, 
\xi^{(1)}_{a,b,\gamma}, \xi^{(2)}_{a,b,\gamma}, 
\xi^{(3)}_{a,b,\gamma}\}$ for small $a$, $b$, and $\gamma$, by 
extending the bases above. Notice first that the vectors in 
the basis for $\gamma=0$ satisfy 
\[ 
  \mathcal S\xi^{(0)}_{a,b,0} \,=\, - \xi^{(0)}_{a,b,0},\quad \mathcal 
  S\xi^{(1)}_{a,b,0} \,=\, - \xi^{(1)}_{a,b,0},\quad \mathcal 
  S\xi^{(2)}_{a,b,0} \,=\, \xi^{(2)}_{a,b,0},\quad \mathcal 
  S\xi^{(3)}_{a,b,0} \,=\, \xi^{(3)}_{a,b,0}~, 
\] 
where $\mathcal S$ is the reversibility operator (\ref{e:s}). 
Since for $\gamma\not=0$ we have $\mathcal S\mathcal A\abg = 
-\mathcal A_{a,b,-\gamma}\mathcal S$, the vectors in the basis can 
be taken such that 
\[ 
  \mathcal S\xi^{(0)}_{a,b,\gamma} \,=\, - \xi^{(0)}_{a,b,-\gamma},\quad 
  \mathcal S\xi^{(1)}_{a,b,\gamma} \,=\, - \xi^{(1)}_{a,b,-\gamma},\quad 
  \mathcal S\xi^{(2)}_{a,b,\gamma} \,=\, \xi^{(2)}_{a,b,-\gamma},\quad 
  \mathcal S\xi^{(3)}_{a,b,\gamma} \,=\, \xi^{(3)}_{a,b,-\gamma}~, 
\] 
and then the matrix $\mathcal M_{a,b,\gamma}$ satisfies 
\[ 
  \widetilde{\mathcal S}\mathcal M_{a,b,\gamma} \,=\, - 
  \mathcal M_{a,b,-\gamma}\widetilde{\mathcal S}~,\quad \hbox{where} 
  \quad \widetilde{\mathcal S} \,=\, \pmatrix{ - \mathbf 1_2 & 
  \mathbf 0_2 \cr \mathbf 0_2 & \mathbf 1_2}~. 
\] 
In addition, since $\mathcal A_{a,b}=\mathcal A_{-a,-b}$, we also 
have $\mathcal M_{a,b,\gamma} = \mathcal M_{-a,-b,\gamma}$. 
Together with the results for $a=b=0$ and $\gamma=0$ we 
conclude that 
\begin{eqnarray*} 
  \mathcal M_{a,b,\gamma}  \,=\, \pmatrix{4\rmi\gamma 
  \left( \mathbf D_2 + \rmO(a^2+b^2)\right) & 
  \mathbf M_2(a,b) - 4\gamma^2 \left( \mathbf 1_2 + \rmO(a^2+b^2)\right) \cr 
  4\gamma^2 \left(\mathbf 1_2+ \rmO(a^2+b^2)\right) & 4\rmi\gamma 
  \left( \mathbf D_2+ \rmO(a^2+b^2)\right)}~. 
\end{eqnarray*} 
 
To end the proof we show that the four eigenvalues of 
this matrix are purely imaginary. The structure of the matrix 
$\mathcal M\abg$ implies that its characteristic polynomial is of 
the form 
 \[ 
  \lambda^4 + \rmi \gamma c_3 \lambda^3 + \gamma^2 c_2\lambda^2 + 
  \rmi\gamma^3 c_1\lambda + c_0\gamma^4~, 
 \] 
in which the coefficients $c_j$ depend upon $a$, $b$, and 
$\gamma$. The four roots of this polynomial are symmetric with 
respect to the imaginary axis, because the spectrum of 
$A_{a,b,\gamma}$ is symmetric with respect to the imaginary axis, 
so that the coefficients $c_j$ are real functions of $a$, $b$, 
$\gamma$. In addition, the spectral equalities 
\reff{e:symg1}--\reff{e:symg2} imply that $c_j$ are even in $a$, 
$b$, and $\gamma$, and that when replacing $(a,b)$ by $(b,a)$ the 
coefficients $c_0$, $c_2$ do not change, while $c_1$, $c_3$ change 
sign. 
 
We now set $\lambda = \rmi\gamma X$, and obtain the polynomial 
with real coefficients, 
\begin{equation}\label{e:polyP} 
  P(X)\,=\,X^4 +  c_3 X^3 - c_2X^2 - c_1X + c_0~. 
\end{equation} 
At $a=b=0$ the four eigenvalues of $\mathcal M\abg$ are known, 
which then gives 
\[ 
  P\big|_{a=b=0}(X)\,=\,X^4-32(\gamma^2+1)X^2+256(1-2\gamma^2+\gamma^4)~. 
\] 
In addition, using the explicit formulas for the plane waves we 
compute the roots of $P$ when $a=\gamma=0$: 
\[ 
  X^{(1,2)}_b \,=\, -4\pm2\sqrt2 b+5b^2+\rmO(b^3),\quad 
  X^{(3)}_b\,=\,X^{(4)}_b \,=\, 4-7b^2+\rmO(b^3)~. 
\] 
Similarly, when $b=\gamma=0$, we find 
\[ 
  X^{(1,2)}_a \,=\, 4\pm2\sqrt2 a - 5a^2+\rmO(a^3),\quad 
  X^{(3)}_a\,=\,X^{(4)}_a \,=\, -4+7a^2+\rmO(a^3)~. 
\] 
Combining these formulas with the parity properties mentioned above,  
we conclude that 
\begin{eqnarray*} 
  &c_3\,=\,4(b^2-a^2)+\rmO(a^4+b^4+\gamma^4),\quad 
  c_2\,=\,32-88(b^2+a^2)+32\gamma^2+\rmO(a^4+b^4+\gamma^4)~,&\\[1ex] 
  &c_1\,=\,\rmO(a^4+b^4+\gamma^4),\quad 
  c_0\,=\,256-1664(b^2+a^2)-512\gamma^2 +\rmO(a^4+b^4+\gamma^4)~.& 
\end{eqnarray*} 
A direct calculation now gives 
\begin{eqnarray*} 
  P(0) &=& 256 + \rmO(a^2+b^2+\gamma^2)>0~, \\ 
  P(X^{(4)}_b) &=& -512a^2-1024\gamma^2 + \rmO((a^2+\gamma^2) 
  (a^2+b^2+\gamma^2))<0~,\\ 
  P(X^{(4)}_a) &=& -512b^2-1024\gamma^2 + \rmO((b^2+\gamma^2) 
  (a^2+b^2+\gamma^2))<0~, 
\end{eqnarray*} 
for $a$, $b$, and $\gamma$ sufficiently small.  This shows that 
the polynomial $P$ has four real roots, so that the four 
eigenvalues of $\mathcal A_{a,b,\gamma}$ are purely imaginary. 
This concludes the proof. 
\end{Proof}

\subsection{Spectrum for $\gamma$ close to $\frac12$} 
 
In this case, the arguments are similar to the ones for 
$\sigma_2(\mathcal A\abg)$ in Section~\ref{ss:small}, and we shall 
therefore only state the results and omit the proofs. First, we 
have the following result on the eigenvalues of $\mathcal A\ze$. 
 
\begin{Lemma}\label{l:p3} 
There exist positive constants $\gamma_1$, $c$, and $\delta$, such 
that the following hold, for any $\gamma\in[\gamma_1,\frac12]$, 
and $(a,b)$ with $\|(a,b)\|\leq\delta$: 
\begin{enumerate} 
  \item $\rmi\omega\abg^{+,0}\in B(-\rmi;\frac12)$, 
  $\rmi\omega\abg^{-,-1}\in B(\rmi;\frac12)$; 
  \item $\rmi\omega\abg^{-,0},\ \rmi\omega\abg^{+,-1}, \ 
  \rmi\omega\abg^{\pm,n}\notin B(0;\frac52)$, $n\in\Z\setminus 
  \{-1,0\}$; 
  \item $|\rmi\omega\abg^{+,n} - \rmi\omega\abg^{+,p}|\geq c$, 
  $n,p\in \Z\setminus\{0\}$, $p\neq n$, $p\neq -n$; 
  \item $|\rmi\omega\abg^{-,n} - \rmi\omega\abg^{-,p}|\geq c$, 
  $n,p\in \Z\setminus\{-1\}$, $p\neq n$, $p\neq -2-n$; 
  \item $|\rmi\omega\abg^{+,n} - \rmi\omega\abg^{-,p}|\geq c$, 
  $n\in \Z\setminus\{0\}, \  p\in \Z\setminus\{-1\}$. 
\end{enumerate} 
\end{Lemma} 
 
Next, we proceed as in the proof of Proposition~\ref{p:far} and 
obtain: 
 
\begin{Proposition}\label{p:12} 
There exist positive constants $\gamma_1$, $c$, and $\delta$, such 
that for any $\gamma\in[\gamma_1,\frac12]$, and $(a,b)$ with 
$\|(a,b)\|\leq\delta$, the following holds: 
\begin{enumerate} 
\item The spectrum of $\mathcal A\abg$ satisfies 
\[ 
  \sigma(\mathcal A\abg) \,\subset\, \bigcup_{n\in\Z} 
  B(\rmi\omega\abg^{-,n};c)\,\cup\,  \bigcup_{n\in\Z} 
  B(\rmi\omega\abg^{+,n};c)~, 
\] 
in which the balls $B(\rmi\omega\abg^{\pm,n};c)$ are mutually 
disjoints except for pairs 
$(\rmi\omega\abg^{+,n},\rmi\omega\abg^{+,-n})$, 
$n\in\Z\setminus\{0\}$, or 
$(\rmi\omega\abg^{-,n},\rmi\omega\abg^{-,-2-n})$, 
$n\in\Z\setminus\{-1\}$. 
 
\item Inside each ball $B(\rmi\omega\abg^{\pm,n};c)$ the operator 
$\mathcal A\abg$ has at most two eigenvalues, which are purely  
imaginary.  
\end{enumerate} 
\end{Proposition}

\section{The focusing NLS equation} 
\label{s:foc} 
 
We consider in this section the focusing NLS equation 
\begin{equation}\label{e:fnls} 
  \rmi U_t(x,t) + U_{xx}(x,t) + |U(x,t)|^2U(x,t) \,=\, 0~, 
\end{equation} 
in which $x\in\R$, $t\in\R$, and $U(x,t)\in\C$. This equation 
also possesses a family of small periodic waves of the form  
$U_{a,b}(x) = \rme^{-\rmi t}\rme^{\rmi\ell_{a,b}x}P_{a,b}(k_{a,b}x)$,  
but now 
\begin{eqnarray*} 
  \ell_{a,b} &=& {1 \over 4}(a^2-b^2) + \rmO(a^4 + b^4)~,\\[1ex] 
  k_{a,b} &=& 1 + {3 \over 4}(a^2+b^2)  + \rmO(a^4 + b^4)~,\\[1ex] 
  P_{a,b}(y) &=& a \rme^{-\rmi y} + b \rme^{\rmi y} 
  + \rmO(|ab|(|a| + |b|))~. 
\end{eqnarray*} 
Both the equation and the periodic waves have the same symmetry  
properties as in the defocusing case, so that we can investigate 
the stability of this family of periodic waves in an analogous  
way. 
 
As in Section~\ref{s:os}, we define $p_{a,b}$ and $Q_{a,b}(z)$ by 
\reff{e:qQdef}, and consider solutions of \reff{e:fnls} of the form 
\reff{e:Qgen}. The wave profile $Q_{a,b}(z)$ is then an equilibrium of 
the evolution equation 
\begin{equation}\label{e:fQeq} 
  \rmi Q_t + 4\rmi p_{a,b}k_{a,b} Q_z + 4k_{a,b}^2 Q_{zz} + 
  (1-p_{a,b}^2)Q + |Q|^2Q \,=\, 0~. 
\end{equation} 
For the orbital stability, we use the same functional set-up, the 
same conserved quantities $N(Q)$ and $M(Q)$, and the energy 
\[ 
  \EE(Q) \,=\, \int_0^{2\pi} \Bigl( 2k_{a,b}^2 |Q_z(z)|^2 
  - \frac14|Q(z)|^4\Bigr)\dd z~, 
\] 
in which only the sign of the last term has been changed. Following 
the arguments in Section~\ref{s:os} one can show that the result in 
Theorem~\ref{th:orbit} holds in this case, as well. We only mention 
that the Hessian matrix of the function $d_{a,b}$ has now the 
expression: 
\[ 
  \HH_{a,b} \,\eqdef\, \pmatrix{ 
    \frac{\partial^2 d_{a,b}}{\partial \omega^2} & 
    \frac{\partial^2 d_{a,b}}{\partial \omega \,\partial c} \vspace{2mm}\cr 
    \frac{\partial^2 d_{a,b}}{\partial c \,\partial \omega} & 
    \frac{\partial^2 d_{a,b}}{\partial c^2}} 
  \Bigg|_{(\omega,c) = (0,0)} \,=\, \frac{\pi}3 
  \pmatrix{2 & 1 \cr 1 & -1}(\1 + \rmO(a^2{+}b^2))~, 
\] 
so that it has again one positive 
and one negative eigenvalue. 
 
The analysis is also the same for the spectral stability, 
when we study the spectrum of the linear operator 
\[ 
  \mathcal A_{a,b} Q \,=\, 4 \rmi k_{a,b}^2 Q_{zz} - 4p_{a,b}k_{a,b} 
  Q_z + \rmi (1-p_{a,b}^2)Q + 2\rmi |Q_{a,b}|^2Q + \rmi 
  Q_{a,b}^2 \overline Q~. 
\] 
However, in this case the result is different: the periodic waves are 
spectrally unstable. While we do not attempt a complete description 
of the spectrum, we focus here on the existence of unstable 
eigenvalues. It turns out that unstable eigenvalues arise through the 
unfolding of the quadruple zero eigenvalue of the unperturbed operator 
at $a=b=\gamma=0$. These are the eigenvalues of the $4\times4$ matrix 
$\mathcal M_{a,b,\gamma}$ in Proposition~\ref{p:zero}, which is 
obtained here in the same way. This matrix has the same structure, 
\begin{eqnarray*} 
  \mathcal M_{a,b,\gamma}  \,=\, \pmatrix{4\rmi\gamma 
  \left(\mathbf D_2 + \rmO(a^2+b^2)\right) & \mathbf M_2(a,b)  
  - 4\gamma^2 \left( \mathbf 1_2 + \rmO(a^2+b^2)\right) \cr 
  4\gamma^2 \left(\mathbf 1_2+ \rmO(a^2+b^2)\right) & 4\rmi\gamma 
  \left( \mathbf D_2+ \rmO(a^2+b^2)\right)}~, 
\end{eqnarray*} 
but now 
\[ 
  \mathbf M_2(a,b) \,=\,  \pmatrix{2a^2 & 4ab \cr 4ab & 2b^2} 
  + \rmO(a^4+b^4)~. 
\] 
The eigenvalues of $\mathcal M_{a,b,\gamma}$ are of the form  
$\lambda = \rmi\gamma X$, where $X$ is a root of a polynomial 
of the form \reff{e:polyP}. When $a=b=0$ we obtain 
\[ 
  P\big|_{a=b=0}(X)\,=\,X^4-32(\gamma^2+1)X^2+256(1-2\gamma^2+\gamma^4)~, 
\] 
and using the plane waves we find the roots of $P$ when $a=\gamma=0$: 
\[ 
  X^{(1,2)}_b \,=\, -4\pm \rmi 2\sqrt2 b-5b^2+\rmO(b^3),\quad 
  X^{(3)}_b\,=\,X^{(4)}_b \,=\, 4+7b^2+\rmO(b^3)~, 
\] 
and when $b=\gamma=0$: 
\[ 
  X^{(1,2)}_a \,=\, 4\pm \rmi 2\sqrt2 a + 5a^2+\rmO(a^3),\quad 
  X^{(3)}_a\,=\,X^{(4)}_a \,=\, -4-7a^2+\rmO(a^3)~. 
\] 
Then we find the expansions for the coefficients 
\begin{eqnarray*} 
  &c_3\,=\,-4(b^2-a^2)+\rmO(a^4+b^4+\gamma^4),\quad 
  c_2\,=\,32+88(b^2+a^2)+32\gamma^2+\rmO(a^4+b^4+\gamma^4),&\\[1ex] 
  &c_1\,=\,\rmO(a^4+b^4+\gamma^4),\quad 
  c_0\,=\,256+1664(b^2+a^2)-512\gamma^2 +\rmO(a^4+b^4+\gamma^4),& 
\end{eqnarray*} 
which give 
\begin{eqnarray*} 
  P(0) &=& 256 + \rmO(a^2+b^2+\gamma^2)>0~,\\ 
  P(X^{(4)}_b) &=& 512a^2-1024\gamma^2 + \rmO((a^2+\gamma^2) 
  (a^2+b^2+\gamma^2))~,\\ 
  P(X^{(4)}_a) &=& 512b^2-1024\gamma^2 + \rmO((b^2+\gamma^2) 
  (a^2+b^2+\gamma^2))~. 
\end{eqnarray*} 
This {suggests} that the polynomial $P$ has complex roots provided 
$\gamma$ is small compared to $a$ and $b$. In order to prove this,  
we consider the polynomial $P$ when $\gamma=0$ and show that it  
has at least two complex roots. We assume that $b\geq a\geq 0$,  
without loss of generality. Since $P(X) = (X-4)^2(X+4)^2 +  
\rmO(a^2+b^2)$, this polynomial is positive outside two 
$\rmO(b^{1/2})$-neighborhoods of $4$ and $-4$, when $a$ and $b$ 
are sufficiently small. Inside each of these neighborhoods, $P$ 
has at most two real roots. A direct computation gives 
\[ 
  P(-4+Y) \,=\, 512b^2+(512b^2+896a^2)Y 
  +(64-40b^2-136a^2)Y^2-(16-4a^2+4b^2)Y^3+Y^4+\rmO(a^4+b^4)~, 
\] 
from which we conclude that $P$ is positive inside any 
$\rmO(b^{1/2})$-neighborhood of $-4$, for $\|(a,b)\|$ sufficiently 
small. Summarizing, $P$ has at most two real roots, and we 
conclude that the operator $\mathcal A\abg$ has at least one pair 
of eigenvalues off the imaginary axis, for $\gamma$ sufficiently 
small. In view of the symmetry with respect to the imaginary axis 
of the spectrum of $\mathcal A\abg$, one of these eigenvalues has 
positive real part. This proves that the small periodic waves are 
spectrally unstable in this case.

\section*{Appendix: Spectrum of $H_{a,b}$}
 
\renewcommand{\theequation}{A.\arabic{equation}} 
\renewcommand{\theLemma}{A.\arabic{Lemma}} 
\renewcommand{\theRemark}{A.\arabic{Remark}} 

In this Appendix we discuss the spectrum of the linear 
self-adjoint operator $H_{a,b}$ defined in \reff{e:Habdef}. As in 
Section~\ref{s:ss}, we decompose the elements of our function 
space into real and imaginary parts, and work with the matrix 
operator 
\[ 
  H_{a,b} \,=\, \pmatrix{-4k_{a,b}^2\partial_{zz} + 
  (p_{a,b}^2-1)+ 3 R_{a,b}^2 + I_{a,b}^2 & 4p_{a,b}k_{a,b}\partial_z 
  + 2R_{a,b} I_{a,b}\cr - 4p_{a,b}k_{a,b}\partial_z + 2R_{a,b}I_{a,b}& 
  -4k_{a,b}^2\partial_{zz} + (p_{a,b}^2-1) + 
  R_{a,b}^2 + 3I_{a,b}^2}~, 
\] 
where $Q_{a,b} = R_{a,b} + \rmi I_{a,b}$. We prove the following 
result: 
 
\begin{Proposition}\label{p:specH} 
There exists a positive constant $\varepsilon_0$ such that for all 
$(a,b)$ with $\|(a,b)\|\leq\varepsilon_0$, the spectrum of the 
matrix operator $H_{a,b}$ in the Hilbert space of $2\pi$-periodic 
functions $L^2_\per([0,2\pi],\C^2)$ verifies 
\[ 
  \sigma(H_{a,b}) \,=\, \{0,\ \lambda^{(2)}_{a,b},\ \lambda^{(3)}_{a,b} 
  \}\;\cup\; \sigma_1(H_{a,b}) ,\quad \sigma_1(H_{a,b}) \,\subset\, 
  [6,+\infty)~, 
\] 
where $0$ is a double eigenvalue and $\lambda^{(j)}_{a,b}$, $j=2,3$,  
are simple real eigenvalues with 
\[ 
  \lambda^{(2)}_{0,0} \,=\, \lambda^{(3)}_{0,0}\quad \mbox{and} \quad 
  \lambda^{(2)}_{a,b}\,<\, \lambda^{(3)}_{a,b},\ \mbox{for all}\ 
  (a,b)\neq 0. 
\] 
\end{Proposition} 
 
\begin{Proof} Notice first that the parity properties with respect 
to $(a,b)$ of the quantities $k_{a,b}$, $p_{a,b}$, and $Q_{a,b}$  
imply that $\sigma(H_{a,b}) = \sigma(H_{-a,b}) = \sigma(H_{a,-b})$,  
and that $H_{a,b}$ commutes with the symmetry $\mathcal S$ introduced in 
\reff{e:s}. 
 
When $a=b=0$ the operator $H_{a,b}$ reduces to the operator $H_0$ 
in the proof of  Lemma~\ref{th:positive} with spectrum 
$\sigma(H_0) = \{4n(n\pm1),\ n\in\Z\}$, for which $0$ is a quadruple 
eigenvalue and the other eigenvalues are all positive 
and greater or equal to $8$. Then, a standard perturbation 
argument shows that the spectrum of $H_{a,b}$ decomposes as 
\[ 
  \sigma(H_{a,b}) \,=\, \{\lambda^{(0)}_{a,b},\ \lambda^{(1)}_{a,b},\ 
  \lambda^{(2)}_{a,b},\ \lambda^{(3)}_{a,b} \}\;\cup\; 
  \sigma_1(H_{a,b}) ,\quad \hbox{where } \sigma_1(H_{a,b}) 
  \,\subset\, [6,+\infty)~, 
\] 
for $(a,b)$ sufficiently small. The four eigenvalues 
$\lambda^{(j)}_{a,b}$ are the continuation for small $(a,b)$ of 
the quadruple zero eigenvalue of $H_0$. 
 
In order to locate these four eigenvalues, we proceed as in the 
proof of Proposition~\ref{p:zero}: we construct an appropriate 
basis 
$\{\xi^{(0)}_{a,b},\xi^{(1)}_{a,b},\xi^{(2)}_{a,b},\xi^{(3)}_{a,b}\}$ 
for the associated four-dimensional eigenspace, compute the 
$4\times4$-matrix $\mathcal M_{a,b}$ representing the action of 
$H_{a,b}$ on this basis, and finally show that the eigenvalues of 
this matrix have the desired property. When $a=b=0$ we choose again 
the basis 
\[ 
  \xi^{(0)}_{0,0} \,=\, \pmatrix{\sin z \cr \cos z}~,\quad  
  \xi^{(1)}_{0,0} \,=\, \pmatrix{0 \cr 1}~,\quad  
  \xi^{(2)}_{0,0} \,=\, \pmatrix{\cos z \cr -\sin z}~,\quad  
  \xi^{(3)}_{0,0} \,=\, \pmatrix{1 \cr 0}~. 
\] 
For $(a,b)\not=0$, the fact that $H_{a,b}$ commutes with the 
symmetry $\mathcal S$ allows us to choose the vectors in the basis 
such that 
\[ 
  \mathcal S\xi^{(0)}_{a,b} \,=\, - \xi^{(0)}_{a,b},\quad \mathcal 
  S\xi^{(1)}_{a,b} \,=\, - \xi^{(1)}_{a,b},\quad \mathcal 
  S\xi^{(2)}_{a,b} \,=\, \xi^{(2)}_{a,b},\quad \mathcal S\xi^{(3)}_{a,b} 
  \,=\, \xi^{(3)}_{a,b}~, 
\] 
and to conclude that the matrix $\mathcal M_{a,b}$ is of the form 
\begin{eqnarray*} 
  \mathcal M_{a,b} &=&  \pmatrix{\mathbf A_2(a,b) & 0 \cr 
  0 & \mathbf B_2(a,b)}~, 
\end{eqnarray*} 
where $\mathbf A_2(a,b)$ and $\mathbf B_2(a,b)$ are $2\times 
2$-matrices with coefficients of order $\rmO(a^2+b^2)$. 
 
Next, the two vectors $\partial_z Q_{a,b}$ and $\rmi Q_{a,b}$ in the 
kernel of $\mathcal A_{a,b}$ also belong to the kernel of $H_{a,b}$, 
so that we can take $ \xi^{(j)}_{a,b} = \xi^{(j)}_{a,b,0}$, for 
$j=0,1$, where $\xi^{(j)}_{a,b,0}$ are as in from the proof of 
Proposition~\ref{p:zero}. Then $\mathbf A_2(a,b) = 0$, so that zero is 
a double eigenvalue: $\lambda^{(0)}_{a,b}=\lambda^{(1)}_{a,b}=0$. The 
remaining vectors $\xi^{(j)}_{a,b}$, $j=2,3$, and the matrix $\mathbf 
B_2(a,b)$ are computed from the expansions of $k_{a,b}$, $p_{a,b}$, 
and $Q_{a,b}$. We find $\xi^{(j)}_{a,b} = \xi^{(j)}_{0,0} + 
\rmO(a^2+b^2)$ for $j=2,3$, and 
\begin{eqnarray*} 
  \mathbf B_2(a,b) &=&\frac1{2\pi}\pmatrix{ \langle 
  H_{a,b} \xi^{(2)}_{0,0}, \xi^{(2)}_{0,0} \rangle  & 
  \langle H_{a,b} \xi^{(2)}_{0,0}, \xi^{(3)}_{0,0} \rangle \vspace{1mm}\cr 
  \langle H_{a,b} \xi^{(3)}_{0,0}, \xi^{(2)}_{0,0} \rangle  & \langle 
  H_{a,b} \xi^{(3)}_{0,0}, \xi^{(3)}_{0,0} \rangle} 
  + \rmO(|ab|(a^2+b^2))\\[2ex] 
  &=& \pmatrix{2a^2 & 4ab \cr 4ab & 2b^2} + \rmO(|ab|(a^2+b^2))~. 
\end{eqnarray*} 
Since the spectrum of $H_{a,b}$ is the same for all couples  
$(\pm a,\pm b)$, the determinant of $\mathbf B_2(a,b)$ is even in both  
$a$ and $b$, which together with the above formula gives 
\[ 
  {\rm det}\left(\mathbf B_2(a,b)\right) \,=\, -12a^2b^2 + 
  \rmO(a^2b^2(a^2+b^2))<0~. 
\] 
This shows that $\lambda^{(2)}_{a,b}<0<\lambda^{(3)}_{a,b}$, for 
sufficiently small $\|(a,b)\|$, which concludes the proof. 
\end{Proof} 
 
\begin{Remark} 
We obtain the same result in the focusing case considered in 
Section~\ref{s:foc}, when the operator $H_{a,b}$ is given by  
\[ 
  H_{a,b} \,=\, \EE_{a,b}''(Q_{a,b}) \,=\, -4k_{a,b}^2 \partial_{zz} 
  -4\rmi p_{a,b}k_{a,b}\partial_z - (1{-}p_{a,b}^2) 
  - |Q_{a,b}|^2 - 2 Q_{a,b} \otimes Q_{a,b}~. 
\] 
The only difference in the proof is the 
expression of the matrix $\mathbf B_2(a,b)$ which is now 
\[ 
  \mathbf B_2(a,b) \,=\, \pmatrix{- 2a^2 & - 4ab \cr - 4ab 
  & - 2b^2} + \rmO(|ab|(a^2+b^2))~, 
\] 
but has the same determinant ${\rm det}\left(\mathbf B_2(a,b)\right) = 
-12a^2b^2 + \rmO(a^2b^2(a^2+b^2))<0$. 
\end{Remark}

\end{document}